\documentclass[2023,preprint,12pt]{elsarticle}

\usepackage{amssymb}
\usepackage{amsmath}
\usepackage{soul,color}
\usepackage{amsfonts}

\usepackage{url}
\usepackage{physics}
\usepackage{xcolor}
\usepackage{caption}
\usepackage{subcaption}
\usepackage{enumitem}
\usepackage[noend]{algpseudocode}
\usepackage{algorithm}
\usepackage{tabularx}
\algnewcommand\algorithmicto{\textbf{to}}
\algnewcommand\algorithmicin{\textbf{in}}
\algnewcommand{\algorithmicand}{\textbf{ and }}
\algnewcommand{\algorithmicor}{\textbf{ or }}
\algnewcommand{\algorithmicnot}{\textbf{not }}
\algnewcommand\algorithmicforeach{\textbf{for each}}
\algrenewtext{For}[3]{\algorithmicfor\ #1 $\gets$ #2\ \algorithmicto\ #3\ \algorithmicdo}
\algdef{S}[FOR]{ForEach}[2]{\algorithmicforeach\ #1\ \algorithmicin\ #2\ \algorithmicdo}
\algnewcommand{\OR}{\algorithmicor}
\algnewcommand{\AND}{\algorithmicand}
\algnewcommand{\NOT}{\algorithmicnot}

\newcommand{\R}{\mathbb{R}}

\renewcommand{\L}{\mathcal{L}}
\renewcommand{\b}{\boldsymbol}

\newcommand{\x}{\b{x}}
\newcommand{\w}{\b{w}}

\newcommand{\h}{\emph{h}}

\usepackage{lineno}

\journal{Journal of Computational Science}

\graphicspath{}

\begin{document}

\begin{frontmatter}



	\title{Spatially dependent node regularity in meshless approximation of partial differential equations}

	\author[mps,ijs]{Miha Rot\corref{cor1}}
	\ead{miha.rot@ijs.si}
	\author[ijs]{Mitja Jančič}
	\ead{mitja.jancic@ijs.si}
	\author[ijs]{Gregor Kosec}
	\ead{gregor.kosec@ijs.si}

	\address[mps]{International Postgraduate School Jožef Stefan, Jamova Cesta 39, 1000 Ljubljana, Slovenia}
	\address[ijs]{Institute Jožef Stefan, Parallel and Distributed Systems Laboratory, Jamova cesta 39, 1000 Ljubljana, Slovenia}

	\cortext[cor1]{Corresponding author}

	\begin{abstract}
		In this paper, we address a way to reduce the total computational cost of meshless approximation by reducing the required stencil size through spatially varying computational node regularity. Rather than covering the entire domain with scattered nodes, only regions with geometric details are covered with scattered nodes, while the rest of the domain is discretized with regular nodes. A simpler approximation can be used in regions covered by regular nodes, effectively reducing the required stencil size and computational cost compared to the approximation on scattered nodes where a set of polyharmonic splines is added to ensure convergent behaviour.
		
		This paper is an extended version of conference paper entitled ``Spatially-varying meshless approximation method for enhanced computational efficiency''~\cite{hybrid_conference} presented at ``International Conference on Computational Science (ICCS) 2023''.
		The paper is extended with discussion on development and implementation of a hybrid regular-scattered node positioning algorithm (HyNP). The performance of the proposed HyNP algorithm is analysed in terms of separation distance and maximal empty sphere radius. Furthermore, it is demonstrated that HyNP nodes can be used for solving problems from fluid flow and linear elasticity, both in 2D and 3D, using meshless methods.

		The extension also provides additional analyses of computational efficiency and accuracy of the numerical solution obtained on the spatially-variable regularity of discretization nodes. In particular, different levels of refinement aggressiveness and scattered layer widths are considered to exploit the computational efficiency gains offered by such solution procedure.
	\end{abstract}


	\begin{highlights}
		\item The total computational cost of the meshless approximation is reduced by using a less robust but computationally more efficient meshless setup on regular nodes and a more expensive, stable setup on scattered nodes.

		\item Dimension independent node placing algorithm that covers regions near geometric details with scattered nodes and the rest of the domain with regular nodes and supports variable node density (\h-refinement) is presented.

		\item The proposed hybrid regular-scattered meshless discretization is demonstrated by solving non-linear natural convection, both in 2D and 3D, and a contact problem in 3D.

	\end{highlights}


\end{frontmatter}


\section{Introduction}
Although the meshless methods are formulated without any restrictions regarding the node layouts, it is generally accepted that quasi-uniformly-spaced node sets improve the stability of meshless methods~\cite{wendland2004scattered,liu2009meshfree}. Nevertheless, even with quasi-uniform nodes generated with recently proposed node positioning algorithms~\cite{slak2019generation,shankar2018robust,van2021fast}, a sufficiently large stencil size is required for stable approximation. A stencil with $n = 2\binom{m+d}{m}$ nodes is recommended~\cite{bayona2017role} for the local Radial Basis Function-generated Finite differences (RBF-FD)~\cite{tolstykh2003using} method in a $d$-dimensional domain for approximation order $m$. The performance of RBF-FD method --- with approximation basis consisting of Polyharmonic splines (PHS) and monomial augmentation with up to and including monomials of degree $m$ --- has been demonstrated with scattered nodes on several applications~\cite{slak2019adaptive,fornberg2015primer,zamolo2019solution}. On the other hand, approximation on regular nodes can be performed with considerably smaller stencil ($n=5$ in two-dimensional domain) using only monomial basis~\cite{kosec2018local} or only Radial Basis Function (RBF)~\cite{kosec2008solution}.

A way to reduce the overall computational complexity while maintaining accuracy is therefore to divide the domain into regions where we need scattered nodes to conform the irregular geometry and the rest of the domain that can be covered with regular nodes that enable using approximation with smaller stencils. While it is not necessary to use regular nodes at all, the more of the domain we can discretise with them without compromising the description of the geometry, the better the expected computational performance.

The spatially varying approximation method has already been introduced in the past to address different problems
with different combinations of methods. A hybrid Finite element method (FEM)-meshless method~\cite{elkadmiri2022hybrid} has been proposed to overcome the issues regarding the unstable Neumann boundary conditions in the context of meshless approximation. FEM has been also coupled with meshless method in~\cite{jaskowiec2023coupling} using approximation constraints to solve Poisson's problem, elasticity and thermo-elasticity problems.
Moreover, the authors of~\cite{ding2004simulation,bourantas2018hybrid} proposed a hybrid of Finite Difference Method (FDM) employed on conventional cartesian grid combined with meshless approximation on scattered nodes to solve flows past a circular cylinder and elasticity problems, respectively.
FDM has been also coupled with meshless in the context of geodynamical simulations~\cite{bourantas2015hybrid}, where authors experimented with combination of Eulerian-Lagrangian schemes. These hybrid approaches are well elaborated, provide stable numerical results and are computationally effective, nevertheless, additional implementation-related burden is required on the transition from mesh based discretisation to scattered nodes~\cite{javed2013hybrid}, contrary to the objective of this paper relying solely on the framework of meshless methods, including the generation of hybrid regular-scattered nodes. The overview of discussed hybrid methods is presented in Table~\ref{table:overview}.

\begin{table}[H]
	\begin{tabular}{|l|l|l|}
		\hline
		\textbf{Methods} & \textbf{Problems tackled}                & \textbf{ref}                 \\ \hline
		meshless/FEM     & nonlinear structural problems            & \cite{elkadmiri2022hybrid}   \\
		meshless/FEM     & Poisson's and thermo-elasticity problems & \cite{jaskowiec2023coupling} \\
		meshless/FDM     & flows past a circular cylinder           & \cite{ding2004simulation}    \\
		meshless/FDM     & elasticity problems                      & \cite{bourantas2018hybrid}   \\
		meshless/FDM     & geodynamical simulations                 & \cite{bourantas2015hybrid}   \\
		meshless/FDM     & flow around solid bodies                 & \cite{javed2013hybrid}       \\  \hline
	\end{tabular}
	\caption{Overview of listed hybrid methods.\label{table:overview}}
\end{table}


{In this paper we first propose a dedicated hybrid regular-scattered node positioning algorithm (HyNP) that is capable of handling irregularities in the domain with scattered nodes, while covering the rest of the domain with regular nodes.} The algorithm is dimension independent, i.e.\ the same algorithm can be used to populate $n$-dimensional domains and inherently supports \h-refinement through the spatially dependent node density.  The performance of the proposed algorithm is evaluated in terms of two metrics that are commonly used to evaluate the discretization quality, i.e.\ the distance to the closest neighbours also referred to as separation distance and the largest empty circle radius. Afterwards, the solution procedure based on such spatially-variable node regularity is analysed in terms of computational efficiency and accuracy of the numerical solution. Compared to the original work~\cite{hybrid_conference}, the analyses in this paper are extended to include different levels of \h-refinement aggressiveness and scattered layer widths to further exploit the efficiency gains offered by such solution procedure, as well as an additional test case from liner elasticity, namely the Boussinesq's contact problem.

The paper is organised as follows: In Section~\ref{sec:fill}, the proposed hybrid regular-scattered node positioning algorithm is described, in Section~\ref{sec:method}, the approximation of linear differential operators using meshless methods is briefly presented, in Section~\ref{sec:examples} the numerical examples are given. Finally, conclusions and future work opportunities are presented in Section~\ref{sec:conclusions}.

\section{Hybrid regular-scattered node positioning algorithm (HyNP) algorithm}
\label{sec:fill}

To obtain the hybrid regular-scattered discretization of a $d$ dimensional domain $\Omega$, we propose an extension for the existing variable density scattered node positioning algorithm proposed by Slak and Kosec~\cite{slak2019generation}. The iterative algorithm begins with a given set of seed nodes that are placed in an ``expansion queue''. In each iteration one node is dequeued and ``expanded''. The expansion stands for a procedure, where several candidate nodes are uniformly generated on a sphere and then randomly rotated around the expanded node. Candidates that do not violate the proximity criteria (are too close to the existing nodes) and are within the domain, are accepted as new nodes and added to the expansion queue. The iteration continues as long as there are nodes in the expansion queue.

\begin{figure}
	\centering
	\includegraphics[width=\textwidth]{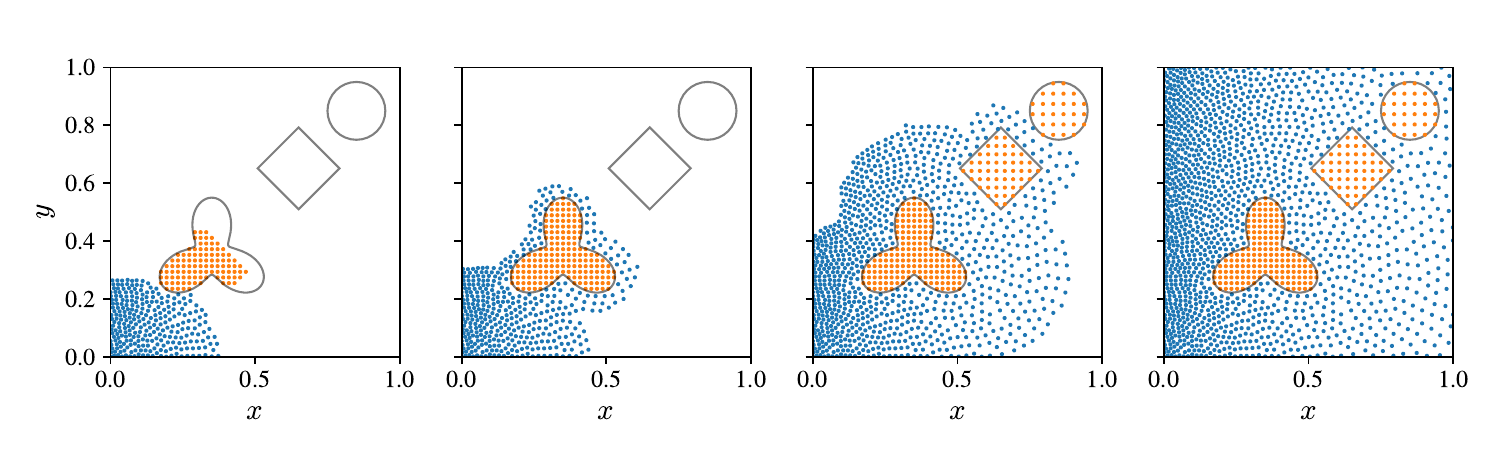}
	\caption{A visualisation of hybrid fill algorithm progression on a domain with variable node density and irregularly shaped areas with regular node positioning.}
	\label{fig:fillProgression}
\end{figure}

In HyNP, we exploit the advancing front nature of the algorithm to find and fill the regular parts of the domain as shown in Figure~\ref{fig:fillProgression}. As soon as the advancing front encounters a regular area defined with an user-defined characteristic function
\begin{equation}
	g\colon\Omega \subset \R^d \to \{0, 1\},
\end{equation}
where $0$ stands for the areas to be populated with scattered nodes and $1$ for the areas to be populated with regular nodes. Regular nodes are placed in a similar advancing front fashion as scattered ones~\cite{slak2019generation}, where the candidate nodes are positioned regularly around the parent node.

The edges of regular area are then used as seed nodes for further progression of the scattered nodes front. Internodal distance in the regular area is determined by the value of the nodal spacing function
\begin{equation}
	h\colon\Omega \subset \R^d \to (0, \infty)
\end{equation}
in the first node that is placed within.

\begin{algorithm}
	\scriptsize
	\caption[Hybrid fill algorithm]{Hybrid fill algorithm.}
	\label{alg:hybridFill}
	\vspace{2pt}
	\textbf{Input:} A $d$ dimensional domain $\Omega$ defined with a characteristic function $\omega\colon \Omega \subseteq \R^d \to \{0, 1\}$. \\
	\textbf{Input:} A nodal spacing function $h\colon\Omega \subset \R^d \to (0, \infty)$. \\
	\textbf{Input:} A characteristic function for regular parts of the domain $g\colon\Omega \subset \R^d \to \{0, 1\}$. \\
	\textbf{Input:} An optional set of boundary and/or seed nodes $X \subseteq \Omega$. \\ 
	\textbf{Output:} A list of nodes in $\Omega$ with regularity based on $g$ and distributed according to spacing function $h$.
	\begin{algorithmic}[1]
		\Function{hybridFill}{$\Omega$, $h$, $g$, $X$}
		\If{$\|X\| = 0$}
		\State \Call{append}{$X$, $\b{p} \in \Omega$} \Comment{Randomly select a seed node if none were provided.}
		\EndIf
		\State $T_p \gets \Call{kdTreeInit}{X}$  \Comment{Initialize spatial search structure on points $X$.}
		\State $T_r \gets \Call{kdTreeInit}{\{\}}$  \Comment{\textbf{Initialize spatial search structure for removable points.}}
		\State $toRemove \gets \{\}$  \Comment{\textbf{List of nodes to remove due to grid conflicts.}}
		\State $i \gets 0$  \Comment{Current node index.}
		\While{$i < |X|$} 	\Comment{Until the queue is not empty.}
		\State $\b{p}_i \gets X[i]$   \Comment{Dequeue current point.}
		\If{$\NOT g(\b{p}_i) \OR g(\b{p}_i) \neq g(\b{p}_{i-1})$} \Comment{\textbf{For scattered and first grid nodes.}}\label{alg:for_scattered_line}
		\State $h_i \gets h(\b{p}_i)$  \Comment{Compute its nodal spacing.}
		\Else
		\State $h_i \gets h_{i-1}$  \Comment{\textbf{Use previous nodal spacing.}}\label{alg:use_previous_line}
		\EndIf
		\ForEach{$\b{c}$}{\Call{candidates}{$\b{p}_i, h_i, g(\b{p}_i)$}} \Comment{\textbf{Generate new candidates.}}
		\label{ln:candidates}
		\If{$\b{c} \in \Omega$}  \Comment{Discard candidates outside the domain.}
		\State $n_p, d_p \gets \Call{kdTreeClosest}{T_p, \b{c}}$ \Comment{Find nearest permanent node index and distance.}
		\If{$d_p \geq h_i $} \Comment{Test that the candidate is far enough.}
		\State $n_r, d_r \gets \Call{kdTreeClosest}{T_r, \b{c}}$ \Comment{\textbf{Find the nearest removable.}}
		\If{$g(\b{c})$} \Comment{\textbf{In regular part of the domain.}}
		\State \Call{prepend}{$X, \b{c}$}  \Comment{\textbf{Enqueue $\b{c}$ as the first element of $X$.}}
		\label{ln:prepend}
		\State \Call{kdTreeInsert}{$T_p, \b{c}$}  \Comment{\textbf{Insert $\b{c}$ into the permanent search structure.}}
		\While{$d_r < h_i$}  \Comment{\textbf{Keep removing the removables while in conflict.}}
		\label{ln:removal}
		\State \Call{append}{$toRemove, n_r$}  \Comment{\textbf{Append to the list of conflicting nodes.}}
		\State \Call{kdTreeRemove}{$T_r, n_r$}  \Comment{\textbf{Remove from the search structure.}}
		\State $n_r, d_r \gets \Call{kdTreeClosest}{T_r, \b{c}}$  \Comment{\textbf{Find the next closest removable.}}
		\EndWhile
		\Else
		\If{$d_p \geq h_i $} \Comment{Test that the candidate is far enough.}
		\State \Call{append}{$X, \b{c}$}  \Comment{Enqueue $\b{c}$ as the last element of $X$.}
		\State \Call{kdTreeInsert}{$T_r, \b{c}$}  \Comment{Insert $\b{c}$ into the removable spatial search structure.}
		\EndIf
		\EndIf
		\EndIf
		\EndIf
		\EndFor
		\State $i \gets i + 1$  \Comment{Move to the next non-expanded node.}
		\EndWhile
		\State \Call{remove}{$X, toRemove$}  \Comment{Remove the conflicting nodes.}
		\State \Return $X$
		\EndFunction
	\end{algorithmic}
\end{algorithm}

The modified parts of the Algorithm~\ref{alg:hybridFill} are highlighted by bold pseudocode comments. The main difference is how the advancing front candidates are generated on line~\ref{ln:candidates}. The algorithm uses $k=15$ randomly placed candidates on a hypersphere around the seed point in the scattered part of the domain and $k = 2d$ candidates $\b c$ at a distance $h$ from $\b p$, computed as given in lines~\ref{alg:for_scattered_line}--\ref{alg:use_previous_line} of Algorithm~\ref{alg:hybridFill} along a standard basis $\widehat {\b e_i}$ that is,
\begin{equation}
	\b{c} = \b{p} \pm h \widehat{\b{e}_i}; \qquad i = {1,...,d}.
\end{equation}
The regular candidate basis could easily be variable throughout the domain, allowing for the regular regions to better match the domain description.
Note that the regular nodes have priority, i.e.\ if a previously placed non-seed scattered node would prevent a grid node to be placed within a regular region, the already accepted scattered node is removed. This ensures that the maximum possible area is covered by one continuous grid of regular nodes and minimizes issues caused by re-entrant grids.

\subsection{Evaluation of discretization quality}
To assess the potential degradation of node quality
due to combining the two different discretization types we assess separation distance metrics of different orders. The first metric is the distance between the $i$-th node and its $j$-th closest neighbour
\begin{equation}
	d_{i, j} = \norm{\b{p}_i - \b{p}_{n(i, j)}},
\end{equation}
where $n(i, j)$ is the index of $j$-th closest neighbour for node $i$. The second metric is a measure of empty space between nodes $s_j$, i.e., the diameters of the largest hyperspheres that can be inscribed in the empty space between generated nodes. The diameters
\begin{equation}
	s_j = 2 \min_i \norm{\b{p}_i - \b{v}_j},
\end{equation}
are then determined by constructing a Voronoi diagram seeded by node positions $\b{p}$ and calculating the distance between the vertex position $\b{v}_j$ and the position of its closest discretisation node $\b{p}_i$ for all vertices in the Voronoi diagram.

\begin{figure}
	\centering
	\includegraphics[width=\textwidth]{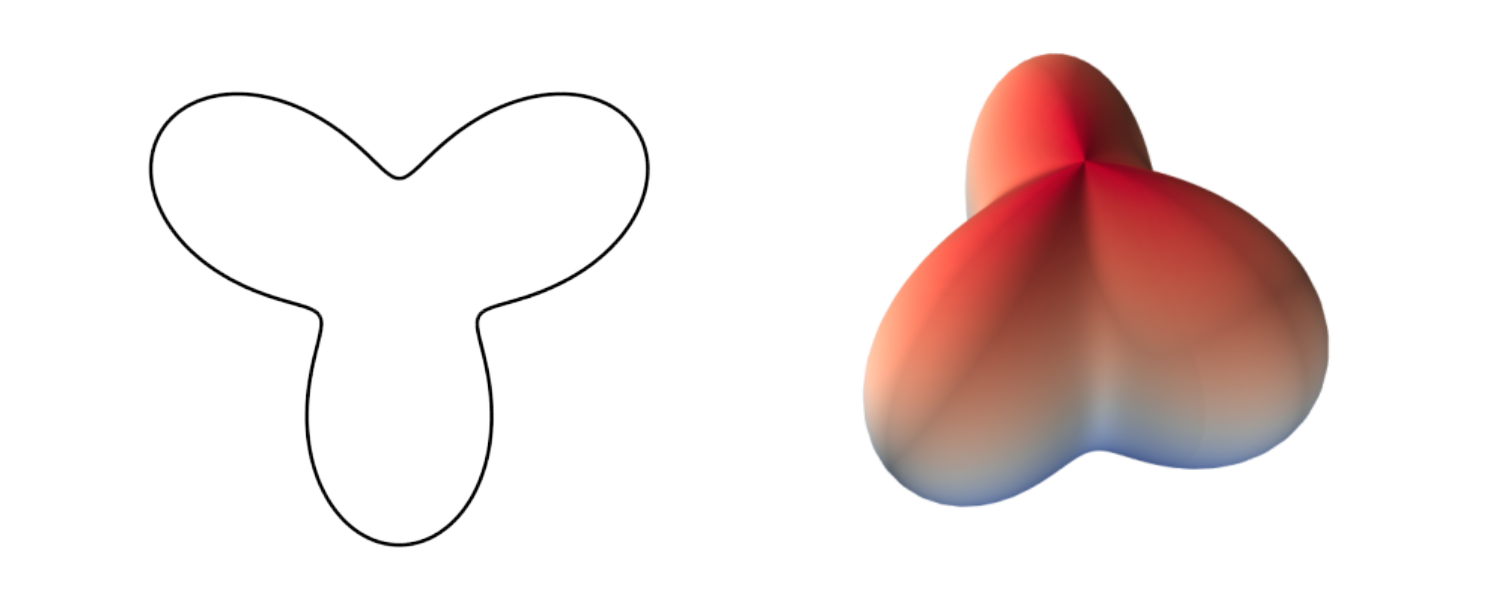}
	\caption{Visualisation of two and three dimensional irregularly shaped clover domains.}
	\label{fig:obstructionShape}
\end{figure}

The two metrics are compared on a square domain with an irregularly shaped grid type characteristic function $g$ based on two and three dimensional clover-like shape shown in Figure~\ref{fig:obstructionShape}. The parametric surface of the boundary $\delta\Theta_{\square}$ of clover-like shape $\Theta_{\square}$, depends on a scaling parameter $l$ and is given by the following expressions
\begin{align}
	r(l, \varphi)      & =\frac{l}{3/2} \Bigg [ 1 - \frac{2}{3}\cos^2 \Bigg(\frac{3}{2}\Big(\varphi - \frac{\pi}{6}\Big )\Bigg ) \Bigg ], \\
	\delta \Theta_{2D} & = r(l, \varphi)\{\cos(\varphi), \sin(\varphi)\} \quad \mathrm{for} \quad \varphi \in \left [ 0, 2 \pi \right ],
	\label{eqn:clover2D}
\end{align}
in 2D and
\begin{align}
	r(l, \varphi, \vartheta) & =\frac{l}{3/2} \Bigg [ 1 - \frac{2}{3}\cos^2 \Bigg(\frac{3}{2}\Big(\varphi - \frac{\pi}{6}\Big )\Bigg ) \frac{\vartheta (\pi -\vartheta)}{3} \Bigg ], &           &                                         \\
	\delta \Theta_{3D}       & = r(l, \varphi, \vartheta)\{\cos(\varphi)\sin(\vartheta), \sin(\varphi)\sin(\vartheta), \cos(\vartheta)\} \quad \mathrm{for}                          & \varphi   & \in \left [ 0, 2 \pi \right ] \nonumber \\
	                         &                                                                                                                                                       & \vartheta & \in \left [ 0, \pi \right ]
	\label{eqn:clover3D}
\end{align}
in 3D with scaling $l = 0.2$. These parametrisations give the boundary between scattered and regular discretization. The parametric definition is practical as it allows for a relatively simple interior check and for the surface $\delta\Theta_{\square}$ to be populated with nodes using a specialized parametric surface node positioning algorithm~\cite{duh_fast_2021}. This clover-like shape will also be used as an irregularly shaped obstacle in further numerical tests.

\begin{figure}
	\centering
	\includegraphics[width=\textwidth]{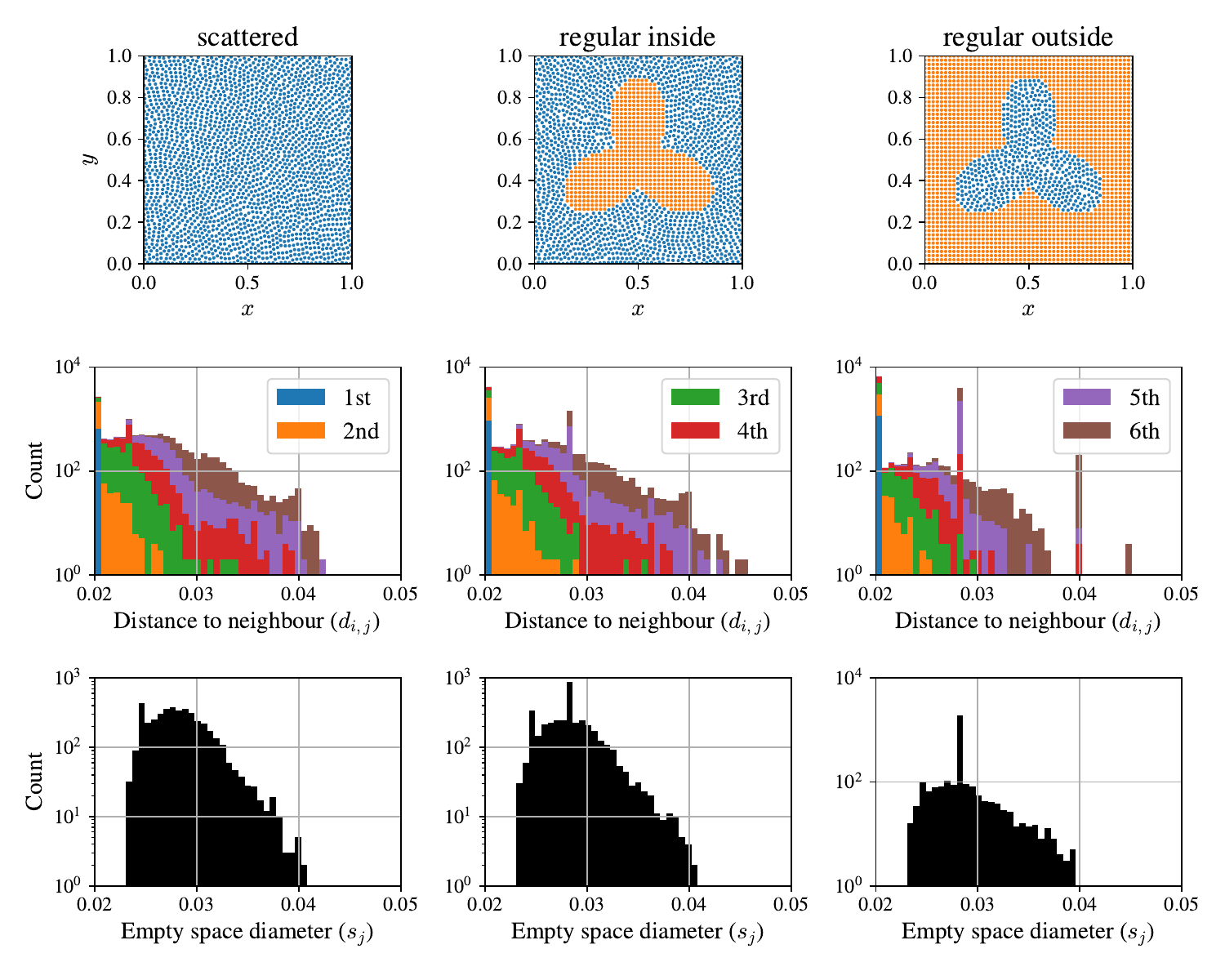}
	\caption{Visualisation of the test domain discretized with different regularity functions $g$ and comparison of fill quality measure distributions between scattered and hybrid fill algorithms on a 2D domain.}
	\label{fig:fillStats2D}
\end{figure}

\begin{figure}
	\centering
	\includegraphics[width=\textwidth]{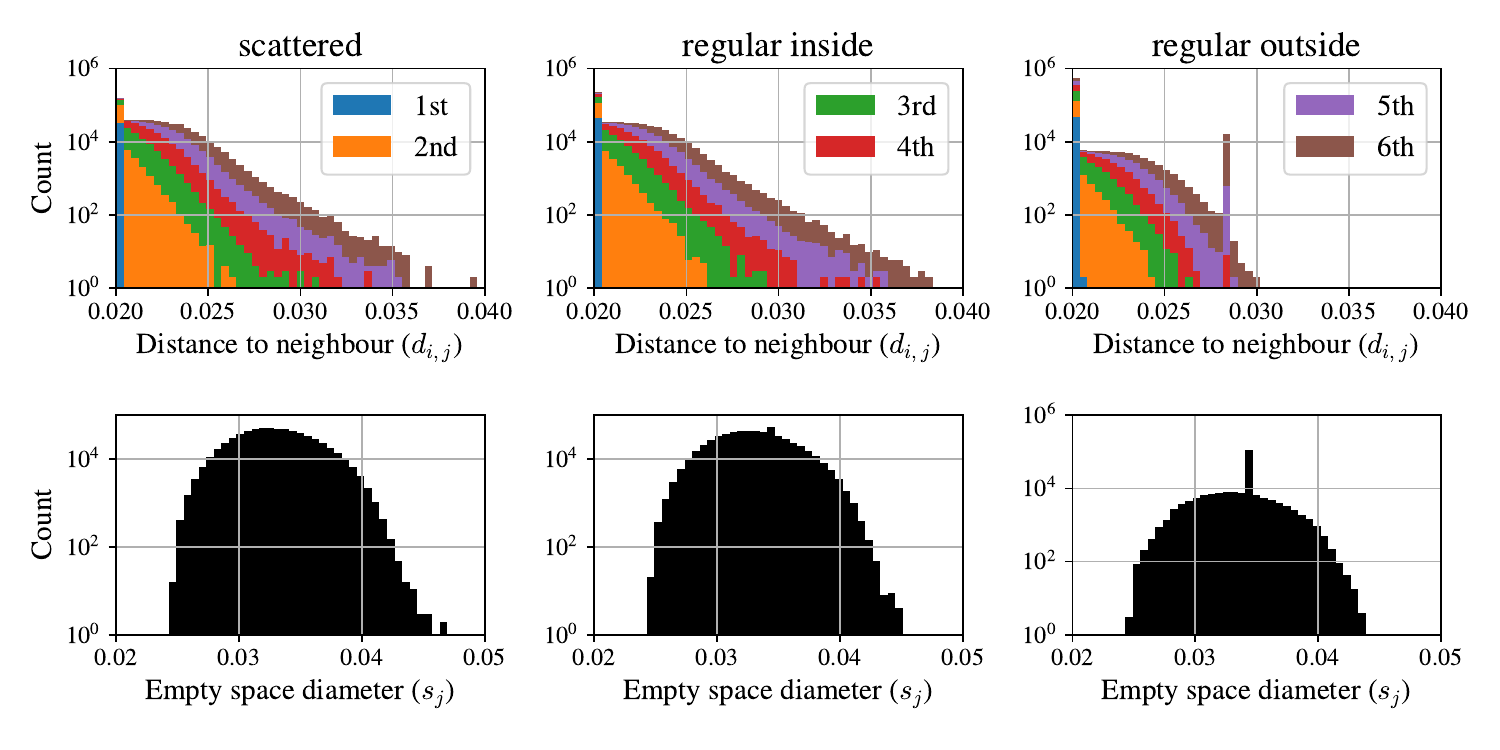}
	\caption{Comparison of fill quality measure distributions between scattered and hybrid fill algorithms on a 3D domain.}
	\label{fig:fillStats3D}
\end{figure}

We first compare three discretizations of a square $\Omega = \left [ 0, 1 \right ] \times \left [ 0, 1 \right ]$ domain shown in the top row of Figure~\ref{fig:fillStats2D}. The first column shows a fully scattered domain, the second column shows a domain filled with regular nodes inside $\Theta$ and scattered outside, while the third column shows the reverse of the second. In all cases the discretization was started with a seed node in the lower left corner and used {a constant} internodal distance $h = 0.02$. The second row shows the distribution of distances to $j$-th neighbour for different discretizations and the third row the radii of the largest possible inscribed circles. We can confirm that the neighbour distance distribution is comparable between the purely scattered and hybrid fill results -- apart from the structural differences that stem from regularity, i.e.\ node clusters at multiples of $h$ and diagonals ($\sqrt{2} h$, $\sqrt{5} h$). More importantly, the empty space distribution also stays the same, confirming that there is no problem in coupling the two node arrangements on the irregular boundary. Additionally we can confirm that there is no discernable difference between the algorithm starting from regular or scattered sections. The analysis is repeated for a cube $\Omega = \left [ 0, 1 \right ] \times \left [ 0, 1 \right ] \times \left [ 0, 1 \right ]$ in Figure~\ref{fig:fillStats3D} with matching conclusions thus confirming the dimensional independence of the algorithm.

\begin{figure}
	\centering
	\includegraphics[width=\textwidth]{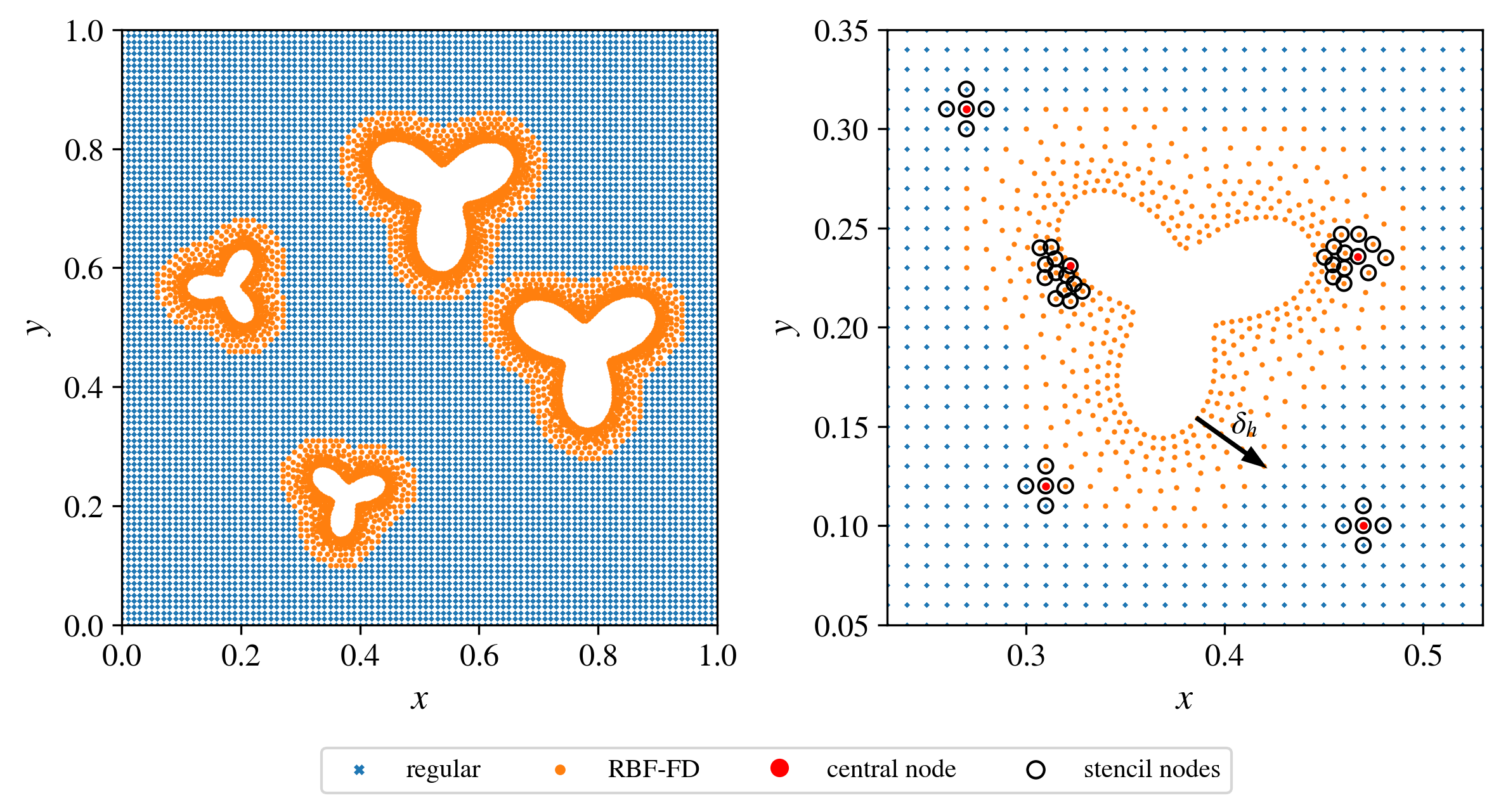}
	\caption{Irregular domain discretization example {\it(left)} and spatial distribution of approximation methods along with corresponding example stencils {\it(right)}.}
	\label{fig:discretization_sample}
\end{figure}

Note that the proposed hybrid discretization is general in the sense that the internodal distance $h$ is by no means limited to a constant, spatially independent, value. Spatially dependent declaration $h(\b p)$ can be used to employ \h-refinement and locally improve the discretization quality where this is necessary. An example of an \h-refined domain discretization with clover-like shaped obstructions is shown in Figure~\ref{fig:discretization_sample}. In this example, the internodal distance linearly increases from $h_s$ on the boundary of the clover to $h_r$ which is equals to the spacing of regularly positioned nodes for a smooth transition. For the demonstration purposes, the width of the scattered node layer $\delta _h$ in Figure~\ref{fig:discretization_sample} has been arbitrarily selected and that the construction of $h(\b p)$ is made through nearest-neighbour search structure, specifically, $k$-d tree.

\section{Numerical approximation of partial differential equations}
\label{sec:method}

With the computational nodes $\x_i\in \Omega$ placed using the HyNP algorithm, the differential operators $\L$ can be locally approximated in point $\x_c$ over a set of $n$ neighbouring nodes (stencil) $\left \{ \x_i \right \}_{i=1}^n = \mathcal{N}$,  using the following expression
\begin{equation}
	\label{eq:ansatz}
	(\L u)(\x _c) \approx \sum_{i=1}^n w_iu(\x _i).
\end{equation}
The approximation~\eqref{eq:ansatz} holds for an arbitrary function $u$ and yet to be determined {vector of} weights $\w$. To determine the weights, the equality of approximation~\eqref{eq:ansatz} is enforced for a chosen set of basis functions. Here we will use two variants
\begin{enumerate}[label=(\roman*)]
	\item The first setup uses a shape parameter free Polyharmonic spline (PHS)
	\begin{equation}
		\varphi(r)=\left\{\begin{array}{ll}
			r^{k}, & k \text { odd } \\
			r^{k} \log r, & k \text { even }
		\end{array}\right. ,
	\end{equation}
	basis augmented with polynomials effectively resulting in a popular \emph{radial basis function-generated finite differences} (RBF-FD) approximation method~\cite{tolstykh2003using}. Such approximation necessarily reproduces polynomials up to the given order (the order of augmenting monomials), in other words, it is exact for augmenting polynomials, i.e. the approximation is of the same order as the polynomial augmentation~\cite{bayona2017role}. This has been discussed, analysed and demonstrated in several recent publications~\cite{janvcivc2021monomial, bayona2017role, bayona2019insight, flyer2016role, fornberg2015primer}, including the recently introduced $hp$-adaptive meshless method, where authors dynamically adjust the order of the method via order of augmenting monomials~\cite{jancic_strong_2023}.

	\item Second setup uses a set Gaussian functions
	\begin{equation}
		g(r)=\exp(-r^2/\sigma^2)
	\end{equation}
	centred at the stencil nodes, where $\sigma$ is shape parameter that has to be determined specifically for each stencil size. We refer to this setup to as a local Radial Basis Function Collocation method (LRBFCM)~\cite{kosec_local_2013, vertnik2013local, kosec2008solution}. In this case, the convergence behaviour is no longer as clear as in the case of augmented approximation. In an experimental study, Ding et al~\cite{ding2005error} showed that the error estimate scales with $O(h/\sigma)$. In~\cite{ding2005error}, authors also discussed the dependence of the error estimate on the stencil size, which was later refined by Bayona et al~\cite{bayona2010}, who found that the method is of second order for stencil sizes between $5$ and $12$. In both papers~\cite{bayona2010, ding2005error} authors experimented with second-order Poisson's PDE in $2$D.
\end{enumerate}

The LRBFCM setup
setup is computationally efficient, but only stable on regular nodes~\cite{kosec2018local,slak2019refinedCauchy}. Unless otherwise specified, in a $d$-dimensional domain the LRBFCM method will be employed using $2d + 1$ Gaussian functions using shape parameter $\sigma=90$. A stencil size equal to the size of the corresponding RBF basis pool is required. For the RBF-FD part, we also resort to the minimal configuration required for 2nd-order operators, i.e., 3rd-order PHS augmented with all monomials up to the 2nd-order ($m=2$) and we can assume that the approximation of partial differential operators used in this work is of a second order~\cite{bayona2010, bayona2017role}. According to the standard recommendations~\cite{bayona2017role}, this requires a stencil size of $n = 2\binom{m+d}{m}$ for a stable approximation.

The weights $\b{w}$ are calculated by imposing equality in Equation~\eqref{eq:ansatz} and solving a system of linear equations $\b{A} \b{w} = \b{b}$ for each computational node
\begin{gather}
	\small
	\begin{split}
		\begin{bmatrix}
			\phi_{1, 1} & \cdots  & \phi_{1, n} \\
			\vdots & \ddots & \vdots \\
			\phi_{n, 1} & \cdots  & \phi_{n, n}
		\end{bmatrix}
		\begin{bmatrix}
			w_1 \\ \vdots \\ w_n
		\end{bmatrix}=
		\begin{bmatrix}
			(\L\phi_{1, 1})|_{\b{p}_1} \\
			\vdots \\
			(\L\phi_{1, n})|_{\b{p}_1}
		\end{bmatrix},
	\end{split}
\end{gather}
where
\begin{equation}
	\phi_{i, j} = \Lambda \left( \frac{\lVert \b{p}_i - \b{p}_j \rVert}{\lVert \b{p}_1 - \b{p}_2 \rVert} \right)
\end{equation}
are radial basis functions written with centralized stencil positions $\b{p}$ and normalized with the distance between the central node position $\b{p}_1$ and the closest stencil position $\b{p}_2$. The basis $\Lambda(r) = g(r)$ for LRBFCM and $\Lambda(r) = \varphi(r)$ for RBF-FD. The system for the latter requires an additional augmentation with $s = \binom{m + d}{d}$ monomials $q$ up to the $m$-th order to ensure positive definitness
\begin{equation}
	\begin{gathered}
		\begin{bmatrix}
			\b{A} & \b{Q} \\
			\b{Q}^T & 0 \\
		\end{bmatrix}
		\begin{bmatrix}
			\b{w} \\ \b{\lambda}
		\end{bmatrix}
		=
		\begin{bmatrix}
			\b{b} \\
			\b{c} \\
		\end{bmatrix},
		\\
		\b{Q} = \begin{bmatrix}
			q_1(\b{p}_1) & \cdots  & q_s(\b{p}_1)\\
			\vdots & \ddots & \vdots \\
			q_1(\b{p}_n) & \cdots & q_s(\b{p}_n)\\
		\end{bmatrix}
		,
		\b{c} = \begin{bmatrix}
			(\L q_1)|_{\b{p}_1} \\
			\vdots \\
			(\L q_s)|_{\b{p}_1}
		\end{bmatrix}.
	\end{gathered}
\end{equation}
The redundant part of the weight vector $\lambda$ is discarded after computation.

It is important to note the difference in required stencil sizes --- 5 vs.~12 nodes in 2D --- that only increases in higher dimensions (7 vs.~30 in 3D). This results both in faster computation of the weights $\w$ (an $\mathcal{O}(N^3)$\footnote{$N_{\mathrm{RBF-FD}} \sim 3 N_{\mathrm{LRBFCM}}$ due to the larger stencil size and the extra PHS in the approximation basis.} operation performed only once for each stencil), in faster evaluation for the $\mathcal{O}(n)$ explicit operator approximation~\eqref{eq:ansatz} performed many times during the explicit time stepping, and in faster solving of the sparse linear systems.

\subsection{Computational stability}

By enforcing the equality of approximation~\eqref{eq:ansatz}, we obtain a linear system $\mathbf{M}\b w = \b \ell$. Solving the system provides us with the approximation weights $\b w$, but the stability of such procedure can be uncertain and is usually estimated via the condition number $\kappa(\mathbf{M}) = \left \| \mathbf{M} \right \| \left \| \mathbf{M}^{-1} \right \|$ of matrix $\mathbf{M}$, where $\left \| \cdot \right \|$ denotes the $L^2$ norm.

\begin{figure}
	\centering
	\includegraphics[width=\textwidth]{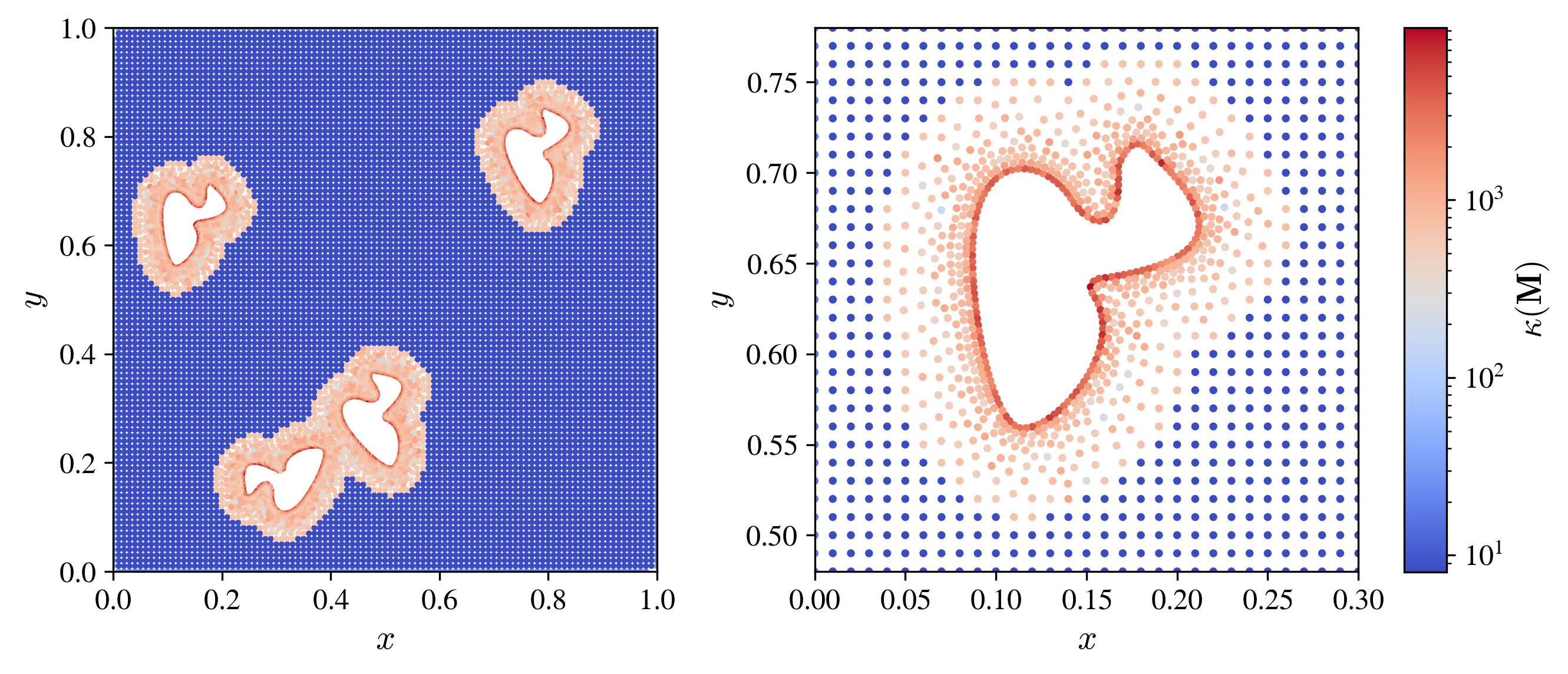}
	\caption{Condition numbers $\kappa(\mathbf{M})$ for the Laplacian operator: entire computational domain {\it(left)} and a zoomed-in section around the irregularly shaped obstacle {\it(right)}.}
	\label{fig:condition_numbers}
\end{figure}

A spatial distribution of condition numbers is shown in Figure~\ref{fig:condition_numbers}. It can be observed that the RBF-FD approximation method generally results in higher condition numbers than the LRBFCM approach. This could be due to the fact that the matrices $\mathbf{M}$ for the RBF-FD part are significantly larger and based on scattered nodes. Nevertheless, it is important to observe that the transition from regular to scattered nodes does not appear to affect the conditionality of the matrices.

\subsection{Implementation details}

We used \texttt{g++ 11.3.0 for Linux} to compile the code with \texttt{-O3 -DNDEBUG} flags on \texttt{Intel(R) Xeon(R) CPU E5520} computer. To improve the timing accuracy we run the otherwise parallel code in a single thread with the CPU frequency fixed at 2.27 GHz, disabled boost functionality and assured CPU affinity using the \texttt{taskset} command.
Post-processing was done using Python 3.10.6 and Jupyter notebooks, also available in the provided git repository\footnote{Source code is available at \url{https://gitlab.com/e62Lab/public/2023_cp_iccs_hybrid_nodes} under tag \textit{v1.3}. }.

\section{Numerical examples}
\label{sec:examples}
\subsection{Natural convection problem}
To objectively assess the advantages of the hybrid discretization method, we first address non-linear natural convection problem that is governed by a system of three PDEs that describe the continuity of mass, the conservation of momentum and the transfer of heat
\begin{align}
	\div \b{v}                                    & = 0, \label{eq:physics1}                                                                            \\
	\pdv{\b{v}}{t} + \b{v} \cdot \grad{\b{v}} & = -\grad p + \frac{1}{\text{Re}}\div(\grad \b{v}) - \b{g} T_\Delta, \label{eq:physics2} \\
	\pdv{T}{t} + \b{v} \cdot \grad{T}             & = \frac{1}{\text{Re}\text{Pr}}\div( \grad T), \label{eq:physics3}
\end{align}
$\b{v}$ is the velocity vector, $p$ the pressure, $T$ the temperature, and $T_\Delta$ the offset from reference temperature. The equations are written in a dimensionless form using Reynolds (Re) and Prandtl (Pr) numbers~\cite{fusegi1991numerical, kosec2008solution} while the results are expressed in terms ob the Rayleigh (Ra) number using the $\mathrm{Ra} = \mathrm{Re}^2 \mathrm{Pr}$ relation. 

The temporal discretization of the governing equations is solved with the explicit Euler time stepping where we first update the velocity using the previous step temperature field in the Boussinesq term~\cite{Tritton1988}. The pressure-velocity coupling is performed using the Chorin's projection method~\cite{chorin1968numerical} under the premise that the pressure term of the Navier-Stokes equation can be treated separately from other forces and used to impose the incompressibility condition.
The time step is a function of internodal spacing $h$, and is defined as $\mathrm{d}t = \frac{h}{4}$ to assure stability.


The problem is solved on different geometries employing (i) LRBFCM, (ii) RBF-FD and (iii) their spatially-varying combination. The performance of each approach is evaluated in terms of accuracy of the numerical solution and execution times.

\subsubsection{The de Vahl Davis problem}
\label{sec:dvd}
First, we solve the standard  de~Vahl~Davis benchmark problem~\cite{de1983natural}. The main purpose of solving this problem is to establish confidence in the presented solution procedure and to shed some light on the behaviour of considered approximation methods, the stability of the solution procedure and finally on the computational efficiency. Furthermore, the de~Vahl~Davis problem was chosen as the basic test case, because the regularity of the domain shape allows us to efficiently discretize it using exclusively scattered or regular nodes and compare the solutions to that obtained with the hybrid regular-scattered discretization.

\begin{figure}
	\centering
	\includegraphics[width=0.55\textwidth]{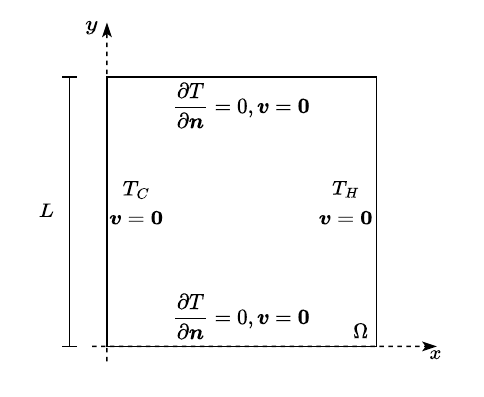}
	\includegraphics[width=0.38\textwidth]{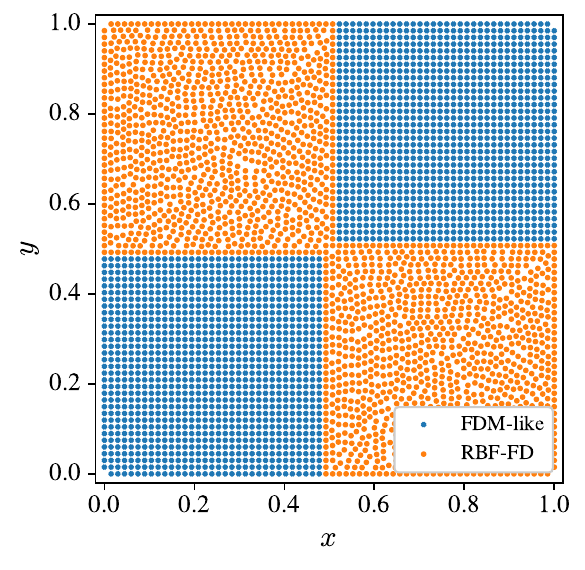}
	\caption{The de Vahl Davis sketch {\it(left)} and example hybrid regular-scattered domain discretization {\it(right)}.}
	\label{fig:dvd_domain}
\end{figure}

For a schematic representation of the problem, see Figure~\ref{fig:dvd_domain} {\it(left)}. The domain is a unit box $\Omega = \left [ 0, 1 \right ] \times \left [ 0, 1 \right ]$, where the left wall is kept at a constant temperature $T_C = -0.5$, while the right wall is kept at a higher constant temperature $T_H = 0.5$. The upper and lower boundaries are insulated, and no-slip boundary condition for velocity is imposed on all walls. Both the velocity and temperature fields are initially set to zero.

To test the performance of the proposed hybrid regular-scattered approximation method, we divide the domain $\Omega$ into quarters, where each quarter is discretized using either scattered or regular nodes -- see Figure~\ref{fig:dvd_domain} {\it(right)} for clarity.

\begin{figure}
	\centering
	\includegraphics[width=\textwidth]{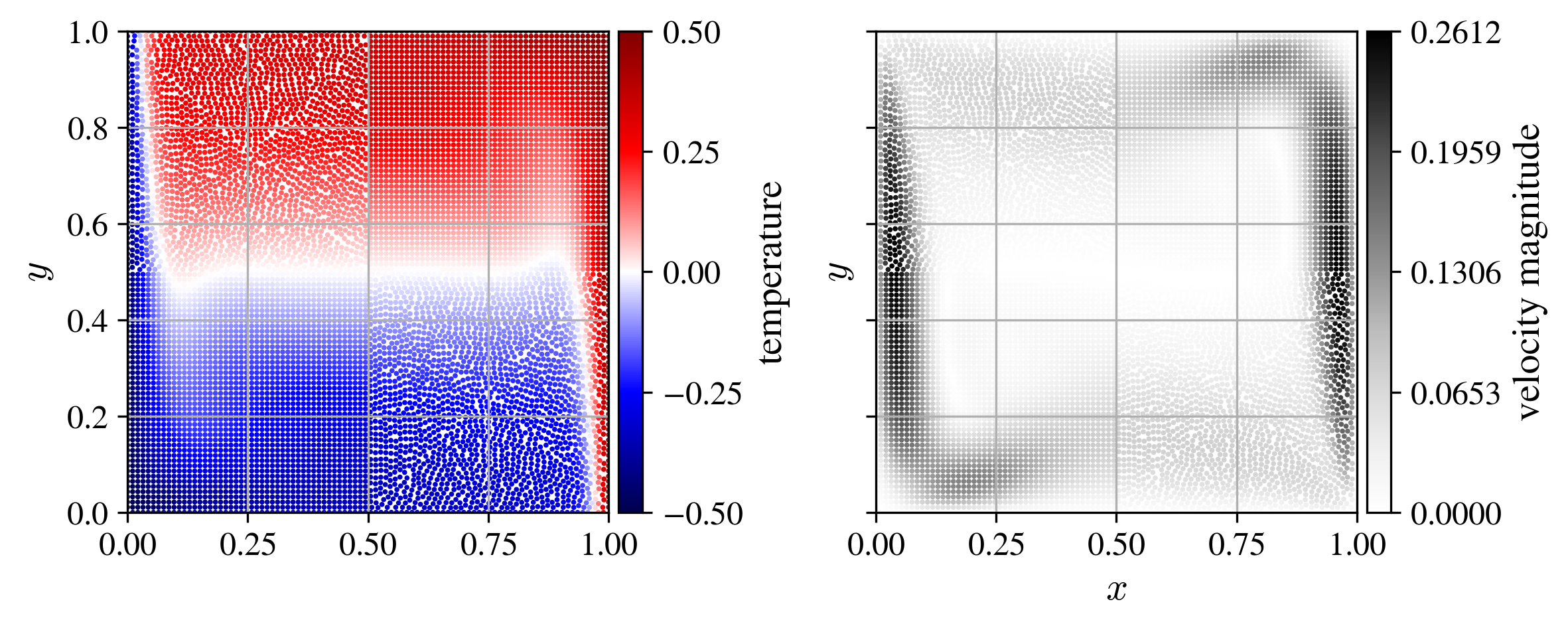}
	\caption{Example solution in the stationary state. Temperature field {\it(left)} and velocity magnitude {\it(right)}.}
	\label{fig:dvd_solution}
\end{figure}

An example solution for $\text{Ra} = 10^6$ and $\text{Pr}=0.71$ at a dimensionless time $t=80$ with approximately $N\approx10\,000$ discretization nodes is shown in Figure~\ref{fig:dvd_solution}.

\begin{figure}
	\centering
	\includegraphics[width=\textwidth]{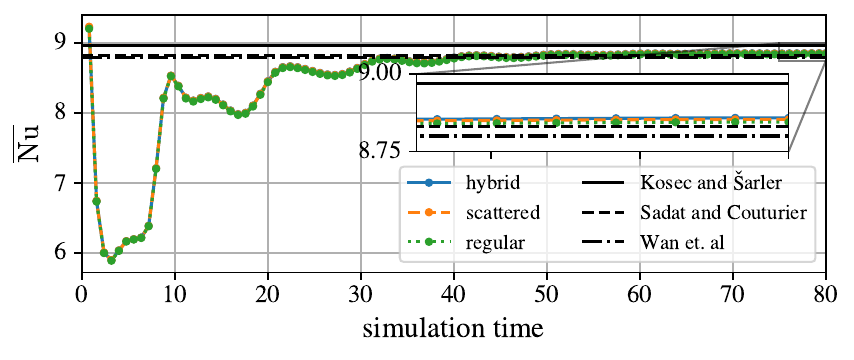}
	\caption{Time evolution of the average Nusselt number along the cold edge calculated with the densest considered discretization. Three reference results Kosec and Šarler~\cite{kosec2008solution}, Sadat and Couturier~\cite{sadat} and Wan et.~al.~\cite{wan} are also added.}
	\label{fig:dvd_nusselt}
\end{figure}

We use the Nusselt number --- the ratio between convective and conductive heat transfer --- to determine when a steady state has been reached and as a convenient scalar value for comparison with reference solutions. In the following analyses, the average Nusselt number ($\overline{\mathrm{Nu}}$) is calculated as the average of the Nusselt values at the cold wall nodes
\begin{equation}
	\mathrm{Nu} = \frac{L}{T_H-T_C}\abs{\pdv{T}{\boldsymbol n}}_{x=0}.
\end{equation}
Its evolution over time is shown in Figure~\ref{fig:dvd_nusselt}. In addition, three reference results are also added to the figure. We are pleased to see that our results are in good agreement with the reference solutions from the literature.

Moreover, Figure~\ref{fig:dvd_nusselt} also shows the time evolution of the average Nusselt number value for cases where the entire domain is discretized using either scattered or regular nodes. We find that all --- hybrid, purely scattered and purely regular domain discretizations --- yield results in good agreement with the references. More importantly, the hybrid method shows significantly shorter computational time than that required by the scattered discretization employing RBF-FD, as can be seen in Table~\ref{tab:nusselt} for the densest considered discretization with $h = 0.00364$.

\begin{table}
	\centering
	\renewcommand{\arraystretch}{1.1}
	\begin{tabular}{cccc}
		Approximation                                    & $\overline{\text{Nu}}$ & execution time [h] & N       \\ \hline \hline
		scattered                                        & 8.854                  & 13.7               & 66\,406 \\
		regular                                          & 8.845                  & 6.2               & 76\,172 \\
		hybrid                                           & 8.856                  & 9.8               & 71\,209 \\ \hline
		Kosec and Šarler (2007)~\cite{kosec2008solution} & 8.97                   & /                  & 10201   \\
		Sadat and Couturier (2000)~\cite{sadat}          & 8.828                  & /                  & 22801   \\
		Wan et.~al.~(2001)~\cite{wan}                    & 8.8                    & /                  & 10201   \\
	\end{tabular}
	\caption{Average Nusselt along the cold edge along with execution times and number of discretization nodes.}
	\label{tab:nusselt}
\end{table}

\begin{figure}
	\centering
	\includegraphics[width=\textwidth]{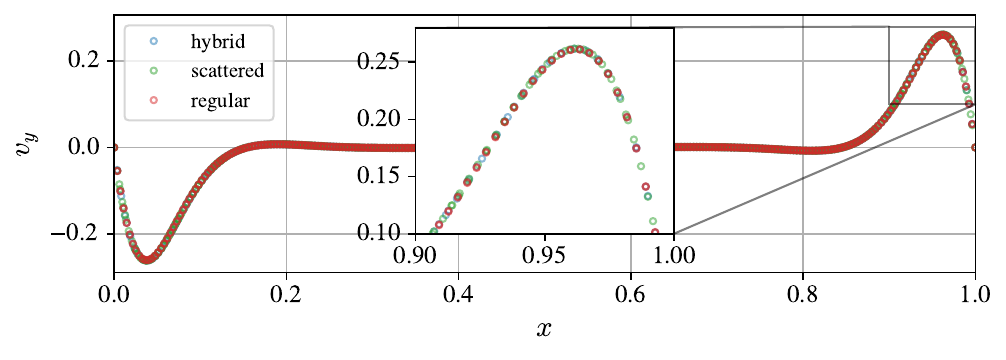}
	\caption{Vertical velocity component values at nodes close to the vertical midpoint of the domain, i.e., $|y - 0.5| \le h$ for purely scattered, purely regular and hybrid discretizations.}
	\label{fig:dvd_sym}
\end{figure}

To further validate the hybrid method, we show in Figure~\ref{fig:dvd_sym} the vertical component of the velocity field across the section $y=0.5$. It is important to observe that the results for the hybrid, scattered and regular approaches overlap, which means that the resulting velocity fields for the three approaches are indeed comparable.

\begin{figure}
	\centering
	\includegraphics[width=\textwidth]{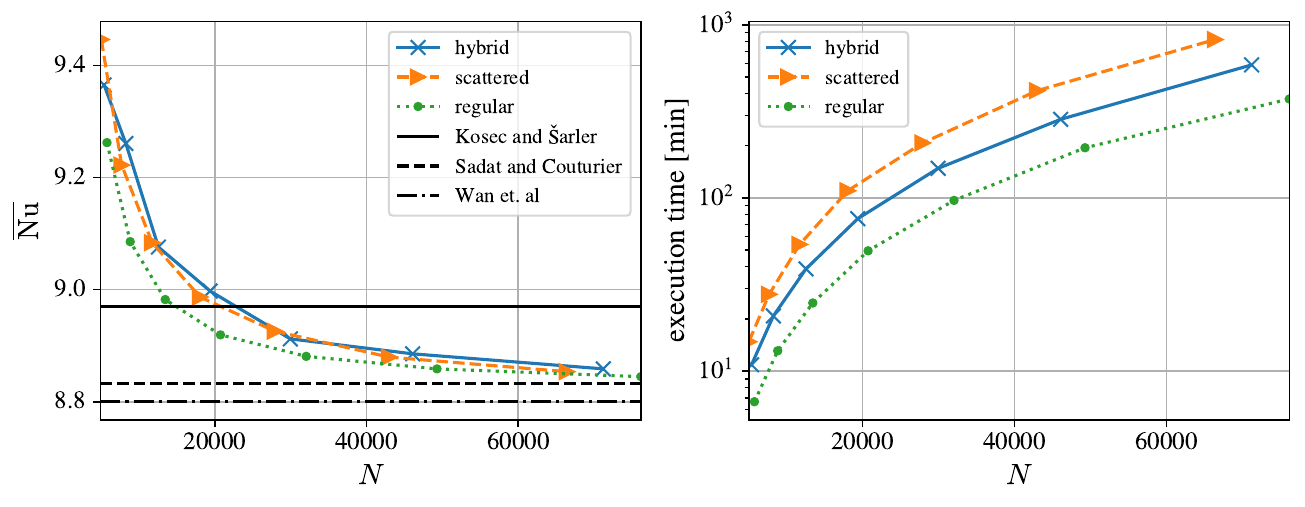}
	\caption{Convergence of average Nusselt number with respect to discretization quality {\it(left)} and corresponding execution times {\it(right)}.}
	\label{fig:dvd_conv}
\end{figure}

As a final remark, we also study the convergence of the average Nusselt number with respect to the number of discretization nodes in Figure~\ref{fig:dvd_conv}, where we confirm that all our discretization strategies converge to a similar value that is consistent with the reference values. Moreover, to evaluate the computational efficiency of the hybrid approach, the execution times are shown on the right. Note that the same values for $h$ were used for all discretization strategies and the difference in the total number of nodes is caused by the lower density of scattered nodes at the same internodal distance.

\subsubsection{The effect of the scattered node layer width $\delta_h$}

To study the effect of the width of the scattered node layer $\delta_h$, we consider two cases. In both cases, the domain from Figure~\ref{fig:dvd_domain} is split into two parts at a distance $h\delta_h$ from the origin in the lower left corner. In the first scenario, the split is horizontal, resulting in scattered nodes below the imaginary split and regular nodes above it. In the second scenario, the split is vertical, resulting in scattered nodes to the left of it and regular nodes to the right of it. In both cases, the domain is discretized with purely regular nodes when $h\delta_h = 0$ and with purely scattered nodes when $h\delta_h = L$.

\begin{figure}
	\centering
	\includegraphics[width=\textwidth]{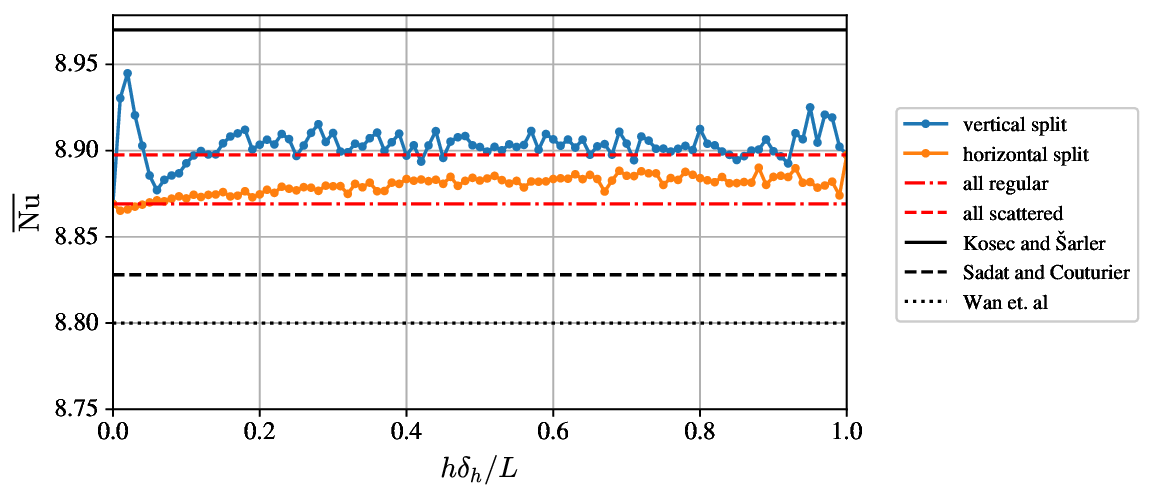}
	\caption{Demonstration of the scattered node layer width ($\delta_h$) effect on the accuracy of the numerical solution.}
	\label{fig:dvd_delta}
\end{figure}

In Figure~\ref{fig:dvd_delta}, we show how the width of the scattered node layer affects the average Nusselt number in stationary state for approximately $40\,000$ discretization nodes. It is clear that even the smallest values of $\delta_h$ yield satisfying results. However, it is interesting to observe that the accuracy is most affected when the boundary between regular and scattered nodes runs across the region with the largest velocity magnitudes, i.e., the first and last couple of vertical split data points in Figure~\ref{fig:dvd_delta}.

\subsubsection{Natural convection on irregularly shaped domains}
\label{sec:irregular}
In the previous section we demonstrated that the hybrid regular-scattered approximation method is computationally more efficient than the pure RBF-FD approximation with only minor differences in the resulting fields. However, to truly exploit the advantages of the hybrid method, irregular domains must be studied. Therefore, in this section, the hybrid regular-scattered approach is employed on an irregularly shaped domain. Let the computational domain $\Omega$ be a difference between the two-dimensional unit box $\Omega = \left [ 0, 1 \right ] \times \left [ 0, 1 \right ]$ and 4 randomly positioned and sized clover-shaped obstacles $\Theta$ defined in Equation~\ref{eqn:clover2D}.

The dynamics of the problem are governed by the same set of equations~(\ref{eq:physics1}-\ref{eq:physics3}) as in the previous section. This time, however, all the boundaries of the box are insulated. The obstacles, on the other hand, are subject to Dirichlet boundary conditions, with half of them at $T_C=-0.5$ and the other half at $T_H=0.5$. The initial temperature is set to $T_{\text{init}} = 0$.

We have chosen such a problem because it allows us to further explore the advantages of the proposed hybrid regular-scattered discretization. Generally speaking, the clover-shaped obstacles within the computational domain represent an arbitrarily complex shape that requires scattered nodes for accurate description, i.e., reduced discretization-related error.

Moreover, by using scattered nodes near the irregularly shaped domain boundaries, we can further improve the local field description in their vicinity by employing a $h$-refined discretization. Specifically, we employ $h$-refinement towards the obstacles with linearly decreasing internodal distance from $h_\mathrm{reg}$ (regular nodes) towards $h_\mathrm{min}$ (irregular boundary) over a distance of $h_\mathrm{reg} \delta_h$. Example discretization is shown in Figure~\ref{fig:discretization_sample} for a scattered node layer width $\delta_h = 4$ and $h_\mathrm{reg}=2h_\mathrm{min}=0.01$, yielding approximately $10\,500$ computational points.

\begin{figure}
	\centering
	\includegraphics[width=\textwidth]{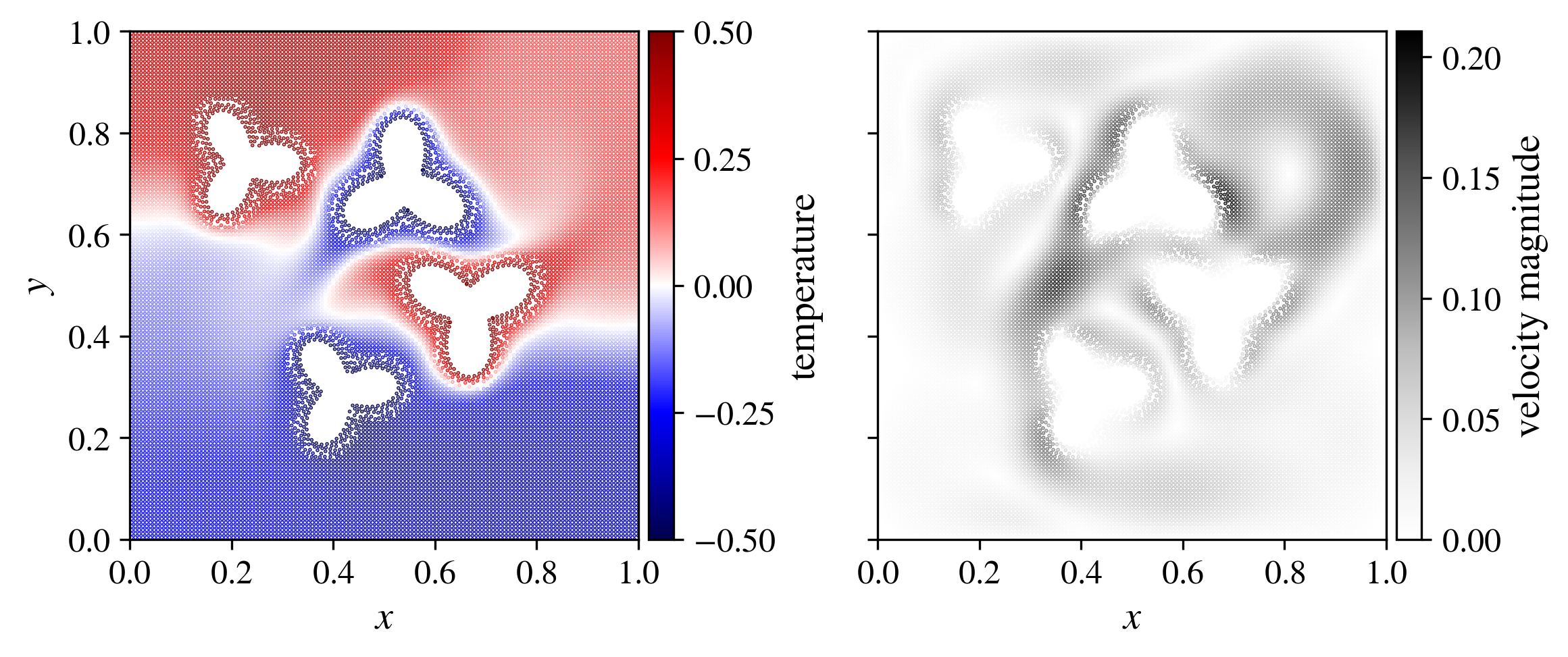}
	\caption{Example solution on irregular domain. Temperature field {\it(left)} and velocity magnitude {\it(right)}.}
	\label{fig:natural_solution}
\end{figure}

\begin{figure}
	\centering
	\includegraphics[width=0.55\textwidth]{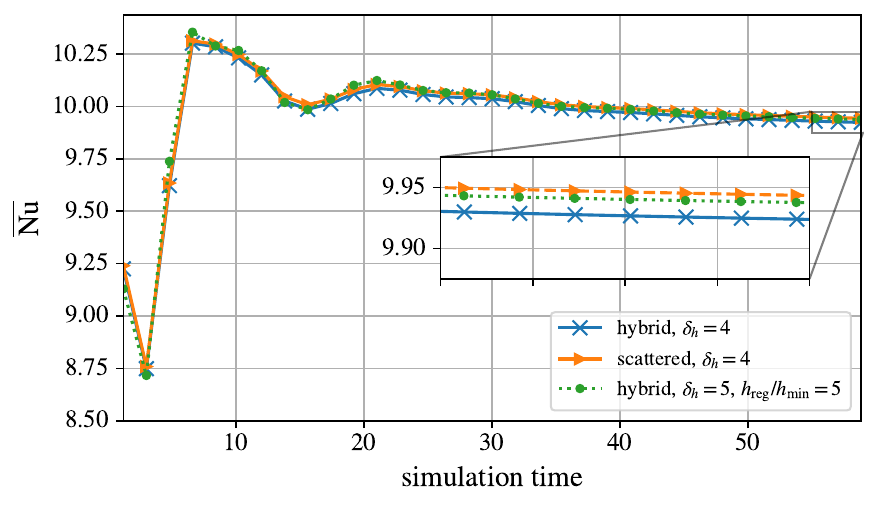}
	\includegraphics[width=0.345\textwidth]{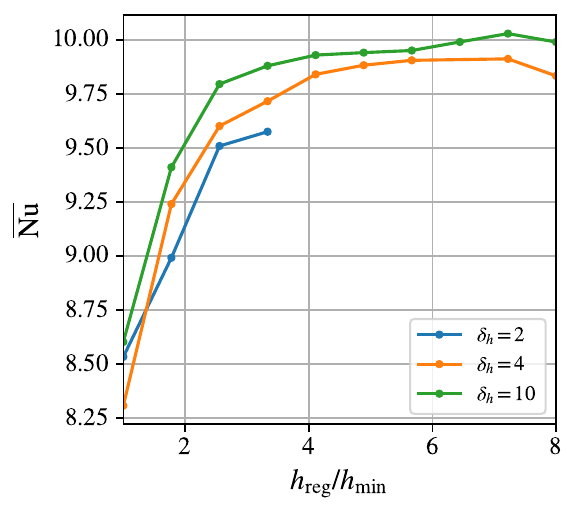}
	\caption{Time evolution of the average Nusselt number calculated on the cold clover-shaped obstacles of an irregularly shaped domain (\it{left}) and average Nusselt number for different combinations of refine parameters (\it{right}).}
	\label{fig:irregular_nusselt}
\end{figure}

Figure~\ref{fig:natural_solution} shows an example solution for an irregularly shaped domain. The hybrid scattered-regular solution procedure was again able to obtain a reasonable numerical solution. Figure~\ref{fig:irregular_nusselt} {\it(left)} shows the average Nusselt number along the cold clover edges where we can observe that a stationary state has been reached. The steady state values for all three considered discretizations match closely. It is perhaps more important to note that the execution times gathered in Table~\ref{tab:nusselt_ireg} show that the uniform density hybrid method effectively reduces the execution time for $\sim$50~\% and that the aggressively refined hybrid discretization for $\sim$90~\%. The purely regular LRBFCM approximation is omitted from the table as it cannot discretise irregular domains.

\begin{figure}
	\centering
	\includegraphics[width=\textwidth]{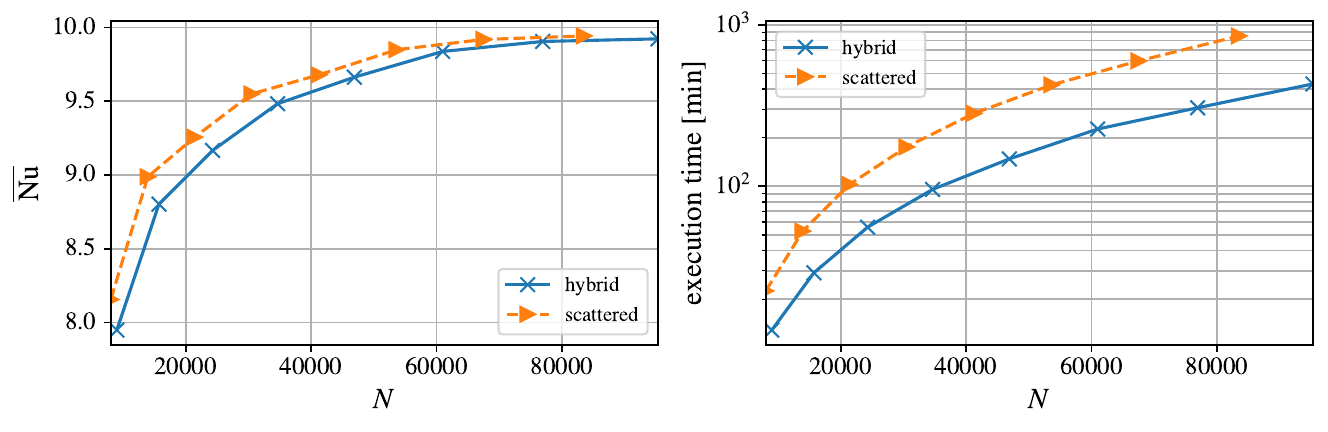}
	\caption{Convergence of average Nusselt number computed on the cold clover-shaped obstacles (left) accompanied with computational times (right).}
	\label{fig:irregular_nusselt_conv}
\end{figure}

The unrefined convergence and computational times are presented in Figure~(\ref{fig:irregular_nusselt_conv}). The results confirm that both the hybrid and the regular discretization converge to a similar Nusselt value and that the hybrid is consistently faster at the same node count while returning a slightly lower value.

\begin{table}[h]
	\centering
	\renewcommand{\arraystretch}{1.1}
	\begin{tabular}{cccc}
		Approximation  & $\overline{\text{Nu}}$ & execution time [h] & N       \\ \hline \hline
		scattered      & 9.942                  & 14.21               & 83\,530 \\
		hybrid         & 9.922                  & 7.17              & 95\,292 \\
		refined hybrid & 9.936                  & 1.52               & 10\,952  \\
	\end{tabular}
	\caption{Average Nusselt along the cold edges of the clover along with execution times. Note that both scattered and hybrid in the table were obtained for $\delta_h = 4$, while refined hybrid was obtained for $\delta_h = 5$ and $h_\mathrm{reg}/h_\mathrm{min} = 5$.}
	\label{tab:nusselt_ireg}
\end{table}

Before continuing with refined discretization we look at interplay between the width of the scattered node layer $\delta_h$ and aggressiveness of refine ratio $\frac{h_\mathrm{reg}}{h_\mathrm{min}}$ shown in Figure~\ref{fig:irregular_nusselt} {\it(right)}. The results confirm our observations from the original work, that the width does not have a significant impact on the result as long as it is wide enough to avoid instability for the selected ratio but there is a slight systematic offset.

\begin{figure}
	\centering
	\includegraphics[width=0.45\textwidth]{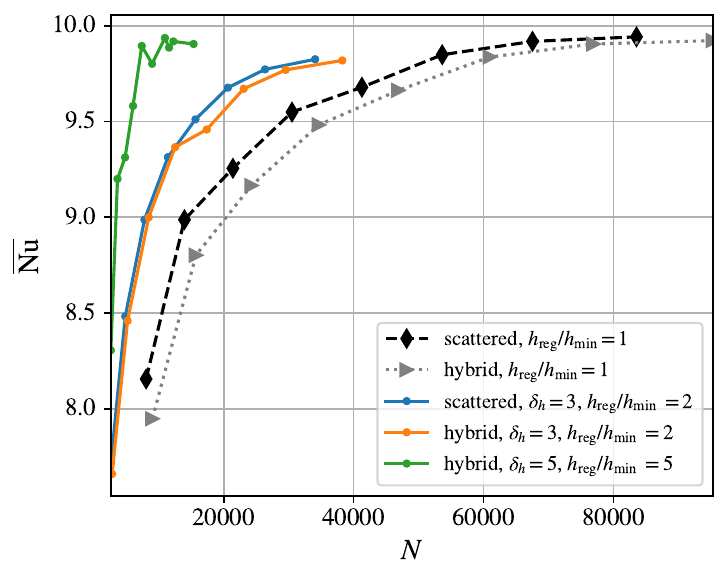}
	\includegraphics[width=0.54\textwidth]{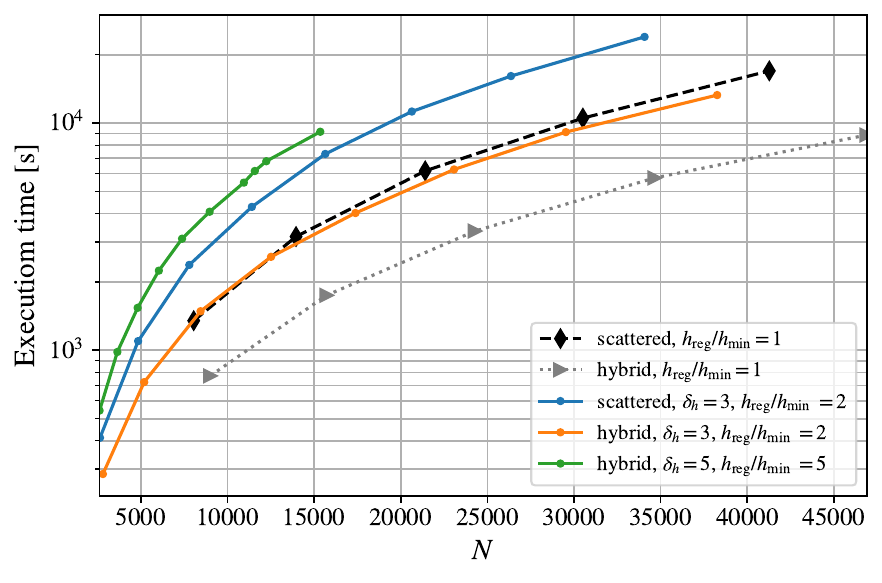}
	\caption{Convergence of average Nusselt number computed on the cold clover-shaped obstacles {\it(left)}  accompanied with computational times {\it(right)}  for refined discretizations.}
	\label{fig:irreg_refine_2d}
\end{figure}

We repeat the convergence study for refined discretizations with results shown in Figure~\ref{fig:irreg_refine_2d}. We chose two hybrid discretizations -- one with less aggressive $\delta_h = 3$, $\frac{h_\mathrm{reg}}{h_\mathrm{min}} = 2$ refine and another larger $\delta_h = 5$, $\frac{h_\mathrm{reg}}{h_\mathrm{min}} = 5$ -- and a scattered discretization with a directly comparable set of parameters. The results on the left graph of Figure~\ref{fig:irreg_refine_2d} show that the refined density solutions are faster to converge to the final Nusselt value, as expected, with the aggressively refined hybrid being significantly faster (also demonstrated in Table~\ref{tab:nusselt_ireg}). The corresponding execution time graph on the right seems more surprising at the first glance due to the refined solutions -- even hybrid ones -- exhibiting comparable or longer times than the uniform scattered discretization but this can be explained by the time-step that is a function of $h_\mathrm{min}$.

\begin{figure}
	\centering
	\includegraphics[width=\textwidth]{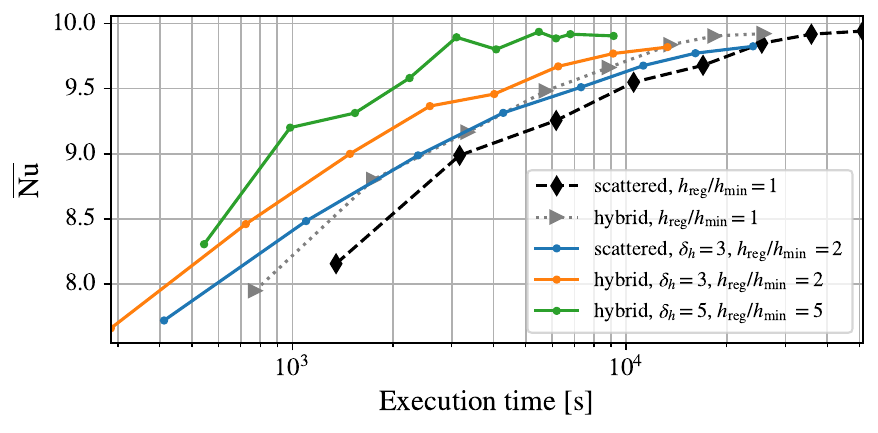}
	\caption{Average Nusselt number calculated on the cold clover shaped obstacles as a function of computational time for a convergence study with a selection of refined and unrefined discretizations.}
	\label{fig:irreg_time_nu_2d}
\end{figure}

The true performance in achieving an accurate solution is easier to determine from a graph of average Nusselt number versus the execution time shown in Figure~\ref{fig:irreg_time_nu_2d}. The results show that the hybrid discretization is slightly faster than the scattered one with comparable refinement and that we can calculate accurate results significantly faster by using a strongly refined hybrid discretization. Surprisingly there is not much difference between weakly refined and unrefined hybrid/scattered approaches.

\subsubsection{Natural convection in three-dimensional domains}
\label{sec:regular}
The de Vahl Davis test is defined on a unit square domain $\Omega = \left [ 0, 1 \right ] \times \left [ 0, 1 \right ]\times \left [ 0, 1 \right ]$.,
where vertical walls are kept at constant temperatures, while
horizontal walls and front/back walls are adiabatic. No-slip velocity
boundary conditions are prescribed on all walls. The dynamics are governed by the same set of Equations (\ref{eq:physics1}-\ref{eq:physics3}) as in the 2D case from Section~\ref{sec:irregular}.

\begin{figure}
	\centering
	\includegraphics[width=\textwidth]{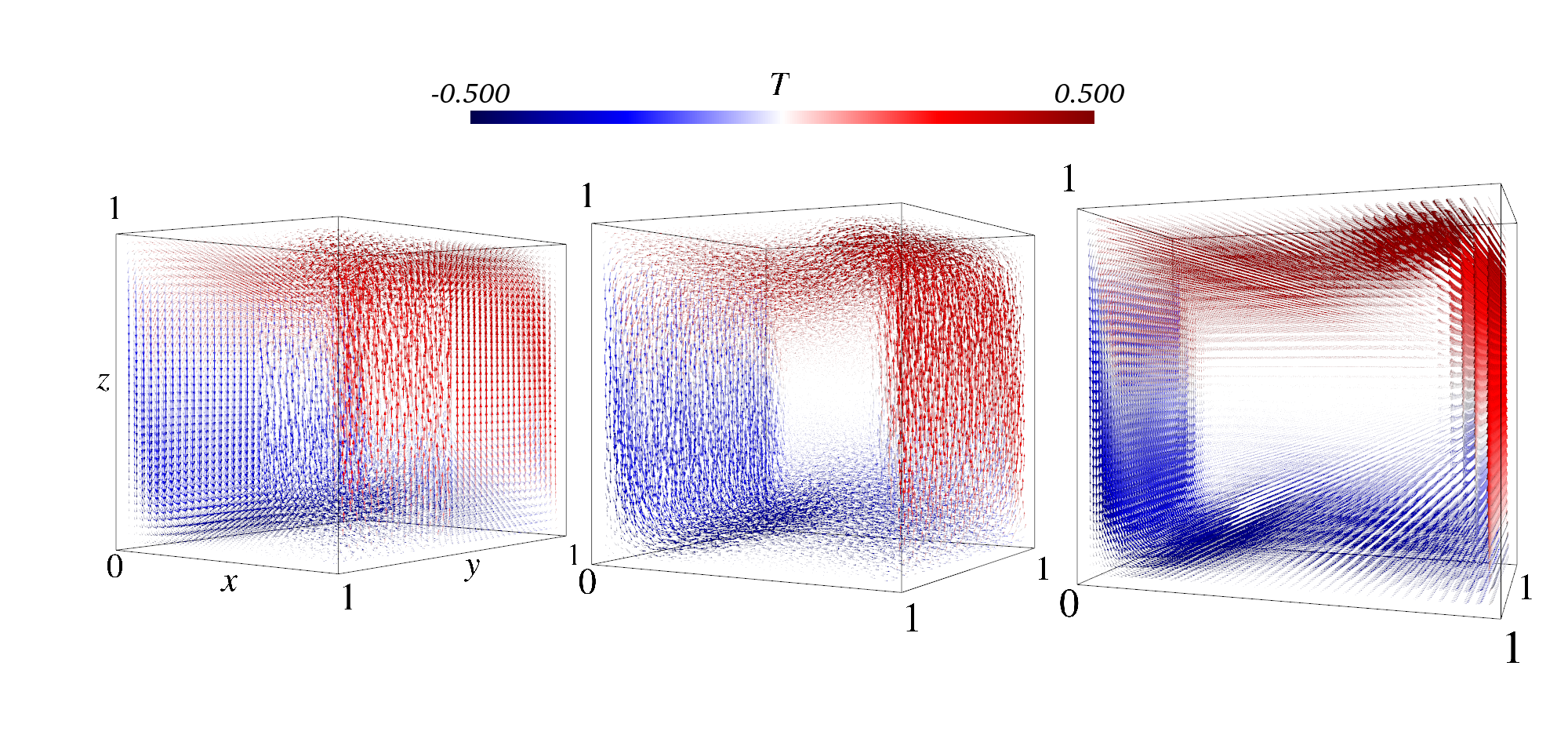}
	\caption{Example of solution of natural convection test using hybrid~(left), scattered~(middle) and regular~(right) nodes in 3D.}
	\label{fig:dvd3d}
\end{figure}

In Figure~\ref{fig:dvd3d} the results for Pr$=0.71$ and Ra$=10^6$ are visualised for all three discretisation variants, namely scattered, regular and hybrid. A more quantitative analysis is presented in Table~\ref{tab:3d_dvd} by comparing characteristic values, i.e. peak positions and values of cross section velocities, with published data.

\begin{table}
	\centering
	\renewcommand{\arraystretch}{1}
	\small
	\begin{tabular}{p{2cm}|cccccc}
		                     Method                      & $v_z({\scriptstyle x_\mathrm{max}, 0.5, 0.5})$ & $x_\mathrm{max}$ & $v_x({\scriptstyle 0.5, 0.5, z_\mathrm{max}})$ & $z_\mathrm{max}$ &   $N$   & $t$ [h] \\ \hline\hline
		                     hybrid                      & 0.2523                                         &      0.960       & 0.0807                                         & 0.133            & 88\,725 & 5.1            \\
		                    regular                      & 0.2322                                         &      0.956       & 0.0801                                         & 0.133            & 96\,800 & 3.1            \\
		                   scattered                     & 0.2592                                         &      0.960       & 0.0817                                         & 0.147            & 81\,218 & 7.9            \\ \hline
		{\footnotesize Slak \& Kosec (2019)~\cite{slak2019generation}}  & 0.2564                                         &      0.961       & 0.0841                                         & 0.143            &    64000    & /              \\
		  {\footnotesize Wang et.\ al (2017)~\cite{wang2017numerical}}   & 0.2556                                         &      0.965       & 0.0816                                         & 0.140            &    125000    & /              \\
		{\footnotesize Fusegi et.\ al (1991)~\cite{fusegi1991numerical}} & 0.2588                                         &      0.966       & 0.0841                                         & 0.144            &    238328    & /              	
	\end{tabular}
	\caption{ Peak positions and values of cross section velocities for 3D natural convection test and number of computational elements. The last column contains the single thread execution times for different discretization strategies.}
	\label{tab:3d_dvd}
\end{table}

As a final demonstrative example of natural convection problems, we employ the proposed hybrid regular-scattered approximation method on a three-dimensional irregular domain, where we add to the domain $\Omega = \left [ 0, 1 \right ] \times \left [ 0, 1 \right ]\times \left [ 0, 1 \right ]$ also 4 randomly positioned and sized clover-like obstacles $\Theta$ defined in Equation~\eqref{eqn:clover3D}.

To improve the quality of the local field description near the irregularly shaped domain boundaries, $h$-refinement is employed with a linearly decreasing internodal distance from $h_r=0.025$ (regular nodes) towards $h_s=0.018$ (clover shapes). The clover-shaped obstacles were set to a constant temperature, two to $T_C = -0.5$ and two to $T_H = 0.5$. The Rayleigh number was set to $10^6$.

\begin{figure}
	\centering
	\includegraphics[width=1.1\textwidth]{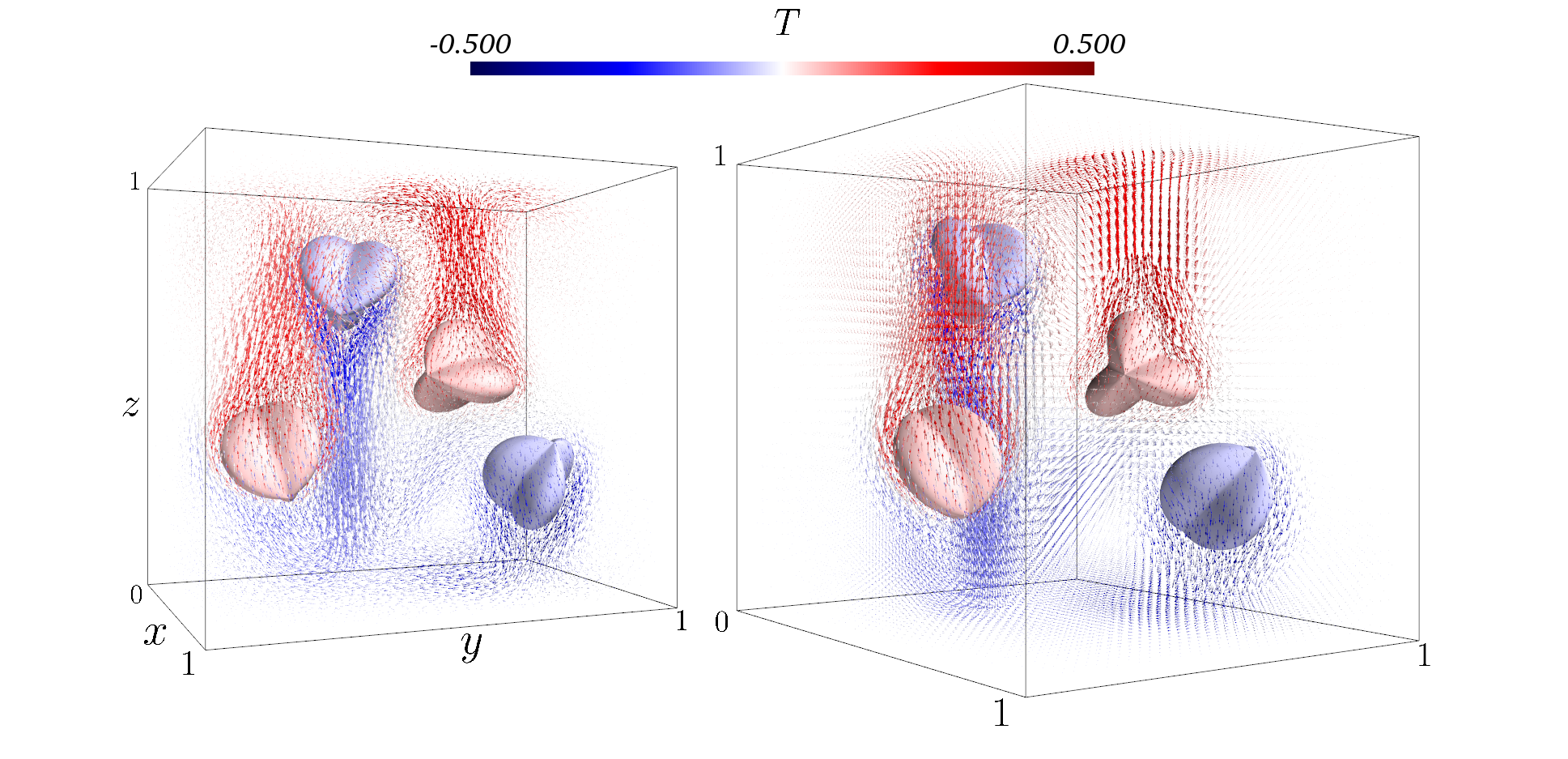}
	\caption{Examples of natural convection in a 3D irregular domain solved using scattered~(left) and hybrid~(right) nodes. The arrows show the velocity in computational nodes and are coloured according to the temperature in that node.}
	\label{fig:3d_solution}
\end{figure}

Although difficult to visualize, an example solution is shown in Figure~\ref{fig:3d_solution}. Using the hybrid regular-scattered domain discretization, the solution procedure was again able to obtain a reasonable numerical solution, i.e.\ the difference in Nusselt number between fully scattered and hybrid approaches is less than $0.5\%$ at the approximately $57\%$ reduction in computation time in favour of the hybrid approach.

\newpage
\subsection{Boussinesq's problem}
\begin{figure}
	\centering
	\includegraphics[width=0.4\textwidth]{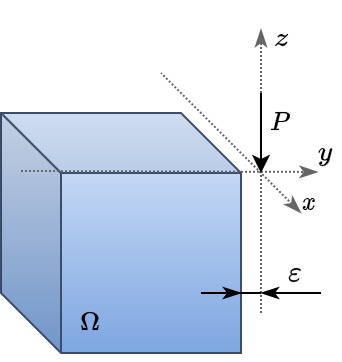}
	\includegraphics[width=0.4\textwidth]{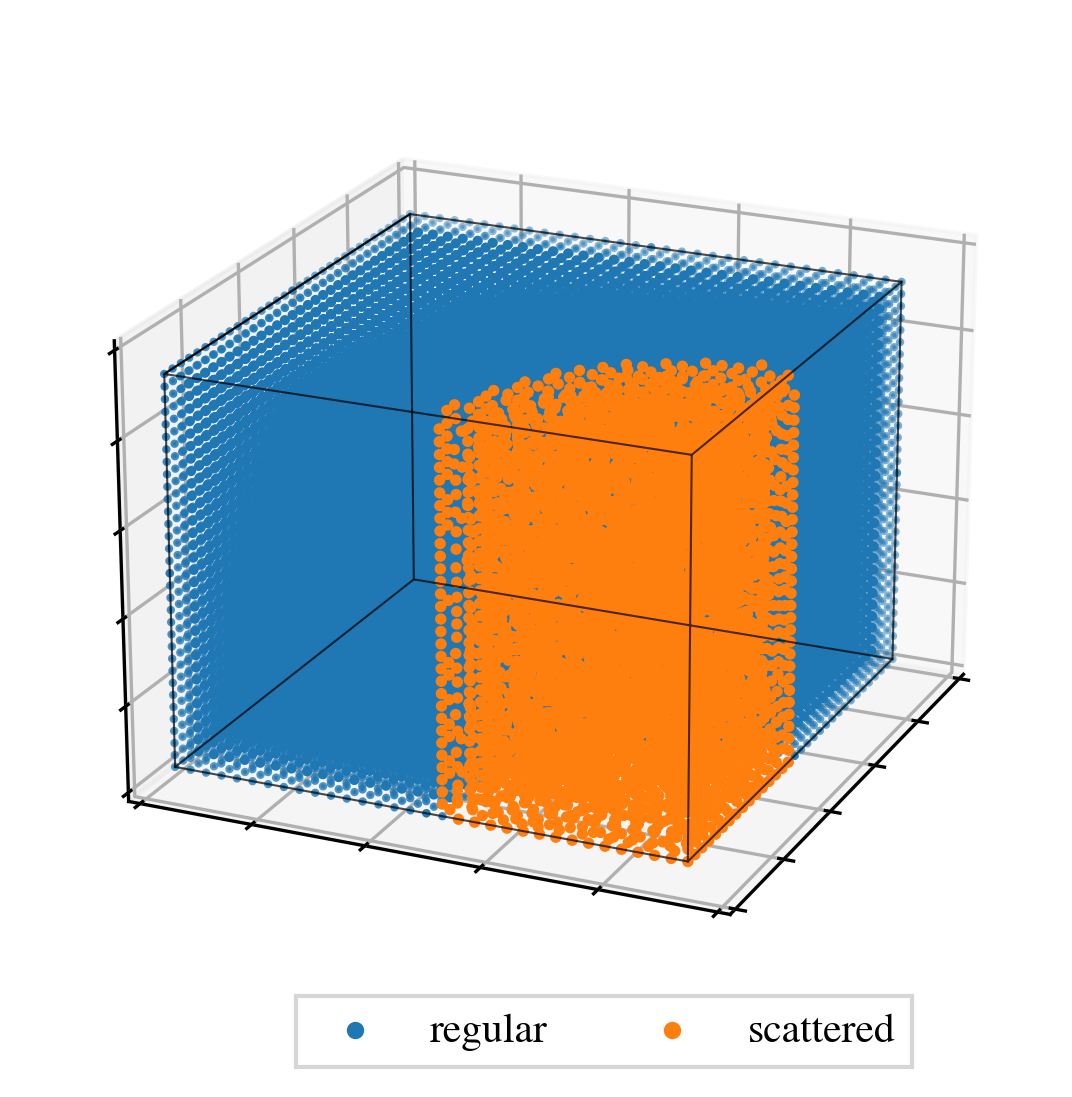}
	\caption{Schematic presentation of Boussinesq's problem {\it(left)} and example spatial distribution of computational node regularity {\it(right)} . For clarity, black lines have been added to mark the domain boundary edges on the right plot.}
	\label{fig:bou_scheme}
\end{figure}

In this section, the proposed solution procedure is demonstrated on an implicit solution to the three-dimensional Boussinesq's problem~\cite{Slaughter_2002}. In this problem, a concentrated normal traction acts on an isotropic half-space, as sketched in Figure~\ref{fig:bou_scheme}.

Analytic solution to the problem is given in cylindrical coordinates $r$, $\theta$ and $z$ as
\begin{align}
	u_r         & = \frac{Pr}{4\pi \mu} \left(\frac{z}{R^3} - \frac{1-2\nu}{R(z+R)} \right), \quad
	u_\theta = 0, \quad
	u_z = \frac{P}{4\pi \mu} \left(\frac{2(1-\nu)}{R} + \frac{z^2}{R^3} \right), \nonumber         \\
	\sigma_{rr} & = \frac{P}{2\pi} \left( \frac{1-2\nu}{R(z+R)} - \frac{3r^2z}{R^5} \right), \quad
	\sigma_{\theta\theta} = \frac{P(1-2\nu)}{2\pi} \left( \frac{z}{R^3} - \frac{1}{R(z+R)} \right),
	\label{eq:3d-problem}                                                                          \\
	\sigma_{zz} & = -\frac{3Pz^3}{2 \pi R^5}, \quad
	\sigma_{rz} = -\frac{3Prz^2}{2 \pi R^5}, \quad
	\sigma_{r\theta} = 0, \quad \sigma_{\theta z} = 0, \nonumber
\end{align}
where $P$ is the magnitude of the concentrated force, $\nu$ is the Poisson's ratio, $\mu$ is the Lam\'e parameter and $R$ is the Eucledian distance to the origin. The solution has a singularity at the origin where the concentrated force is applied to the bulk. This makes the problem difficult to solve and consequently a good candidate for treatment with spatially variable node regularity; allowing us to employ \h-refined scattered nodes towards the singularity and regular nodes elsewhere.

We consider only a part of the domain, i.e.\ $ \varepsilon$ away from the singularity yielding $\Omega$ defined as a box, i.e.\ $\Omega = [-1, -\varepsilon] \times [-1, -\varepsilon] \times [-1, -\varepsilon]$, as schematically shown in Figure~\ref{fig:bou_scheme}. From a numerical point of view, we solve the Navier-Cauchy equation
\begin{equation}
	(\lambda +\mu)\nabla (\nabla \cdot \b u)+\mu\nabla^2\b u = \b f,
	\label{eq:NC}
\end{equation}
using the Lam\'e parameters $\lambda$ and $\mu$ with Dirichlet boundary conditions as given in~\eqref{eq:3d-problem}.

Even though the analytic solution is given in cylindrical coordinate system, the problem is implemented using cartesian coordinates. For the physical parameters of the problem, the values $P=-1$, $E=1$ and $\nu = 0.33$ were used.

To employ the hybrid discretization, a cylinder along the edge with applied force is assumed. Points inside the cylinder are scattered, allowing us to employ \h-refinement towards the critical edge in the continuation of this work, while regular nodes are positioned elsewhere. To determine if a point $\b p$ is inside the cylinder with radius $R_0$, it's perpendicular distance $\left \| (\x_\mathrm{corner} - \b p) \times \widehat{\b e _z}\right \|$ to the edge is computed and compared to $R_0$. Here, $\x_\mathrm{corner}$ is the domain corner closest to the origin, i.e.\ $\x_\mathrm{corner}  =(-\varepsilon, -\varepsilon, -\varepsilon)$ and $\widehat{\b e _z}$ is a unit vector along the $z$ axis. For clarity, spatial distribution of the computational node regularity and of the approximation methods are also shown in Figure~\ref{fig:bou_scheme} {\it(left)}.

Note that the final sparse system was solved using BiCGSTAB with ILUT preconditioner, where the global tolerance was set to $10^{-16}$ with a maximum number of $300$ iterations and drop-tolerance and fill-factor set to $10^{-6}$ and $60$ respectively. Example solution is shown in Figure~\ref{fig:bou_scheme} {\it(left)}  with displacement magnitudes on the left and von Mises stress on the right.

\begin{figure}
	\centering
	\includegraphics[width=\textwidth]{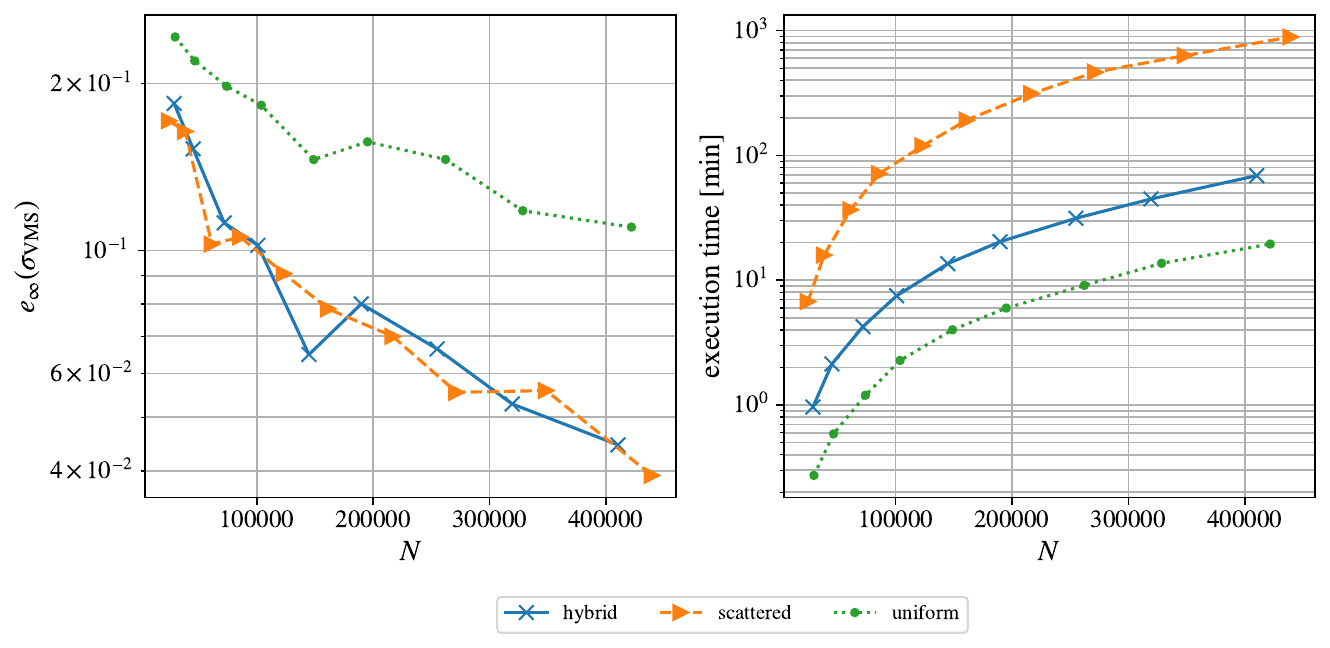}
	\caption{Convergence analysis using regular, hybrid and scattered discretization in case of Boussinesq's problem {\it(left)}  and execution times {\it(right)}.}
	\label{fig:bou_conv_times}
\end{figure}

The convergent behaviour and the computational gains offered by the proposed solution procedure are studied in Figure~\ref{fig:bou_conv_times}. The offset from the origin ($\varepsilon$) was fixed and set to 0.1, the scattered node area was set to $R_0 = 0.4$, and the \h-refinement towards the singularity is avoided by setting \ $h_\mathrm{reg} = h_\mathrm{min}$ for a fair performance comparison of the different domain discretizations approaches. On the left, we show the infinity norm error in terms of von Mises stresses. We observe that using a purely regular discretization, the error is approximately two times larger compared to fully scattered or hybrid discretizations but retains a similar convergence rate.

Figure~\ref{fig:bou_conv_times} shows the wall-clock times required to obtain the numerical solutions. The LRBFCM approach with purely regular domain discretization is clearly the fastest among the three but also the least accurate. Both hybrid and scattered discretizations show approximately an order of magnitude (or more) longer wall-clock times with comparable accuracy of the numerical solution. Note that the parameter $R_0$ strongly impacts the computational time required for the hybrid discretization.

In the following set of analyses, the hybrid and scattered discretizations are employed with a linearly decreasing internodal distance $h(\b p)$, from $h_\mathrm{reg}$ (on regular nodes) to $h_\mathrm{min}$ (on the edge subjected to force $P$) depending on the perpendicular distance from the edge with applied force. Thus, the internodal distance $h$ at point $\b p$ is defined with the following expression:
\begin{equation}
	\label{eq:bou_h}
	h(\b p) = \min \left  \{ h_\mathrm{reg}, h_\mathrm{min} + (h_\mathrm{reg} - h_\mathrm{min}) \frac{\left \| (\x_\mathrm{corner} - \b p) \times \widehat{\ e _z}\right \|}{R_0} \right \}.
\end{equation}

\begin{figure}
	\centering
	\includegraphics[width=\textwidth]{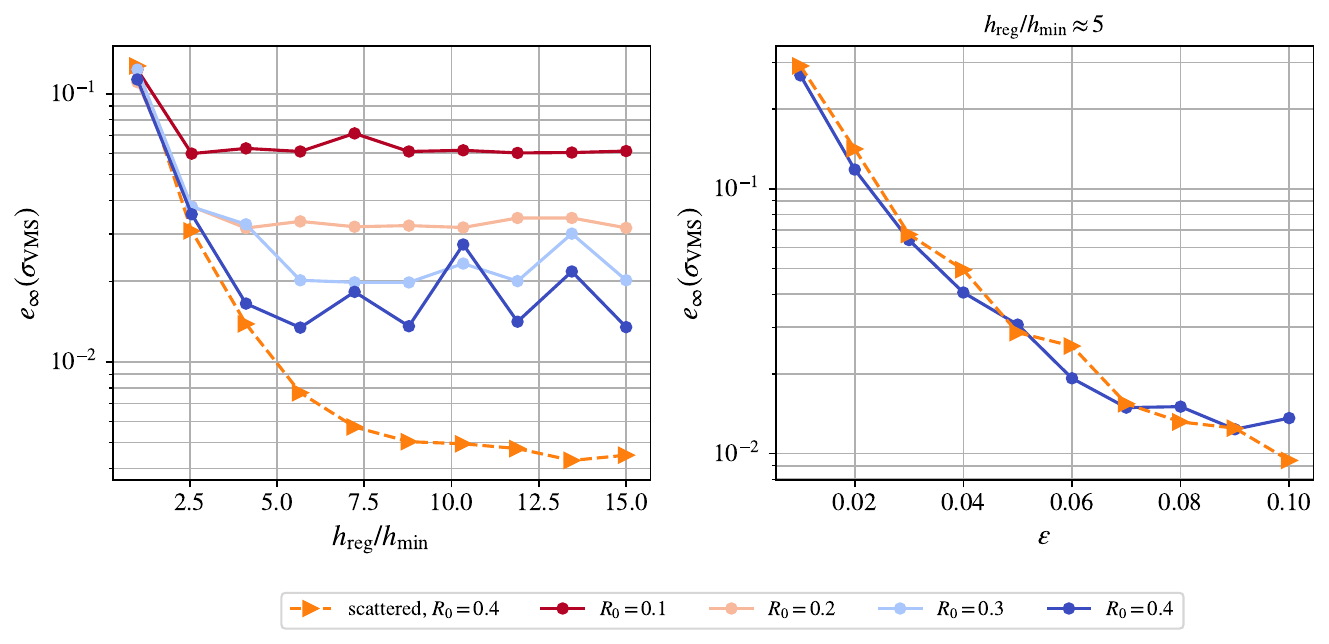}
	\caption{The impact of \h-refine aggressiveness on the maximum error of von Mises stress for scattered and hybrid discretizations with different $R_0$ {\it(left)}  and the maximum von Mises stress as a function of distance $\varepsilon$ between the singularity and the outer edge of the domain {\it(right)} .}
	\label{fig:bou_eps_refine}
\end{figure}

In Figure~\ref{fig:bou_eps_refine} {\it(left)} , we study the behaviour in case of different \h-refinement aggressiveness. The numerical solutions are again evaluated in terms of the infinity norm of the von Mises stress. We show that the \h-refinement towards the edge with applied concentrated force improves the accuracy of the numerical solution by an order of magnitude for both purely scattered and hybrid domain discretization approaches shown for different $R_0$. We did not encounter any stability related issues in the process -- even for the most aggressive \h-refinement used, i.e.\ $h_\mathrm{reg}/h_\mathrm{min}=15$. On the other hand the accuracy reaching a plateau with increasing refinement provides an insight into the trade-offs inherently present in the hybrid method. Once the refined part of the domain is discretized with a sufficiently large density of nodes the accuracy becomes bounded by the regular part of the domain. The error plateau exhibited by the scattered and hybrid methods, otherwise discretized for the same $R_0 = 0.4$, directly reflects the difference in method accuracy observed in Figure~\ref{fig:bou_conv_times}. With decreasing $R_0$, the relatively weak LRBFCM method is used to discretize areas ever closer to the singularity leading to larger errors irrespective of how well the area in the immediate vicinity of the corner is discretized.

In Figure~\ref{fig:bou_eps_refine} {\it(right)} we show a brief study of hybrid and fully scattered discretization approaches with respect to the domain corner distance $\varepsilon$ to the singularity that is present in the origin. As expected, the accuracy of the numerical solutions decreases for smaller distances $\varepsilon$ and improves for larger values. It is worth mentioning, that \h-refinement following Equation~\eqref{eq:bou_h} was used in the process. We used $h_\mathrm{reg}=(1-\varepsilon)/40$ and $h_\mathrm{min} = h_\mathrm{reg}/5$.

\begin{figure}
	\centering
	\includegraphics[width=\textwidth]{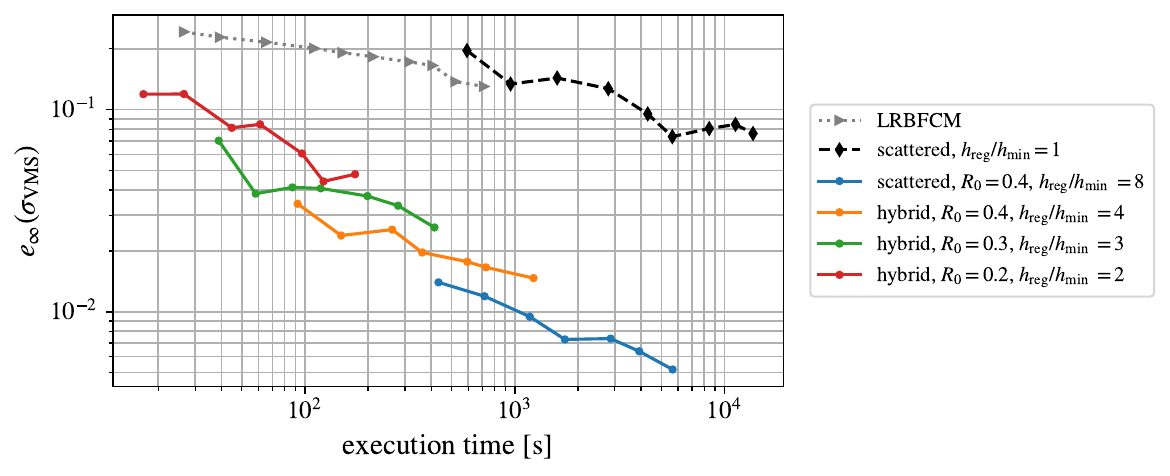}
	\caption{Maximum error of von Mises stress as a function of computational time for a convergence study with a selection of refined and unrefined discretizations.}
	\label{fig:bou_refined_err_time}
\end{figure}

Finally, we can leverage the insights gathered from discretization parameter analysis to repeat the convergence study with a refined scattered and hybrid discretization. We chose 4 configurations to analyse: scattered with $R_0 = 0.4$, $\frac{h_\mathrm{reg}}{h_\mathrm{min}} = 8$ as the most accurate, hybrid with $R_0 = 0.2$, $\frac{h_\mathrm{reg}}{h_\mathrm{min}} = 2$ as the fastest and hybrids with $R_0 = 0.3$, $\frac{h_\mathrm{reg}}{h_\mathrm{min}} = 4$ and $R_0 = 0.4$, $\frac{h_\mathrm{reg}}{h_\mathrm{min}} = 4$ as middle ground. The results for error as a function of computational time are shown in Figure~\ref{fig:bou_refined_err_time}, confirming the expected discretization rankings in both accuracy and execution time. All refined discretizations provide lower error than the regular and scattered constant density solutions at a comparable or lower computational time. In this case a hybrid discretization provides significant time reduction compared to a scattered discretization if we are somewhat willing to compromise on accuracy. Once the discretization is sufficiently refined, the regular part dominates the error and we would be better served by replacing LRBFCM with a more accurate, albeit expensive, approximation for the regular part of the hybrid.

\section{Conclusions}
\label{sec:conclusions}
In this paper we proposed an improvement in terms of computational efficiency for the numerical treatment of problems in which most of the domain can be discretized with regular nodes, while scattered nodes are used only near irregularly-shaped domain boundaries. First, we introduced an algorithm for $n$-dimensional $h$-refined meshless node placement that, based on user input, discretizes different regions either with scattered or with regular nodes. We showed that such an approach does not degrade the quality of generated nodes by means of analysing the separation distance and maximal empty sphere radius width of generated nodes.

The remainder of the paper is dedicated to demonstrating how the proposed hybrid regular-scattered discretization performs in different problems. We combined the regular nodes with LRBFCM, a fast but sensitive method, and scattered nodes with RBF-FD, expensive but robust method. With such a setup we solved the de Vahl Davis natural convection and Boussinesq's contact problems, in 2D and 3D. We showed that the proposed hybrid regular-scattered discretization can significantly contribute to the computational efficiency, while introducing minimal to no cost regarding the accuracy of the numerical solution.

Further analysis is required regarding the selection of an appropriate approximation method for the regular part to avoid the regular part dominating the error as in some of the refined Boussinesq's cases. Additionally, the scattered node layer width and the aggressiveness of $h$-refinement near the irregularly shaped domain boundaries should be investigated, as both affect the computational efficiency and stability of the solution procedure. Future work should also include more difficult problems, such as mixed convection problems and a detailed analysis of possible surface effects, e.g.\ scattering, at the transition layer between the scattered and regular domains.
\subsubsection*{Acknowledgements}
The authors acknowledge the financial support from the Slovenian Research and Innovation Agency (ARIS) research core funding No. P2-0095, Young Researcher programme PR-10468, and research projects No. J2-3048 and No. N2-0275.

Funded by National Science Centre, Poland under the OPUS call in the Weave programme 2021/43/I/ST3/00228.
This research was funded in whole or in part by National Science Centre (2021/43/I/ST3/00228). For the purpose of Open Access,
the author has applied a CC-BY public copyright licence to any Author Accepted Manuscript (AAM) version arising from this submission.

\subsubsection*{Conflict of interest}
The authors declare that they have no conflict of interest. All the co-authors have confirmed to know the submission of the manuscript by the corresponding author.

\bibliographystyle{elsarticle-num}
\bibliography{hybridNodeExtended}
\end{document}